\documentclass{amsart}
\usepackage{amssymb} 
\usepackage{dsfont}
\usepackage{footmisc}
\usepackage{stmaryrd}
\usepackage[usenames,dvipsnames,svgnames,table]{xcolor}
\usepackage[utf8]{inputenc}
\usepackage[T1]{fontenc}
\usepackage{ tipa }
\usepackage{helvet}
\usepackage{color}
\usepackage{mathrsfs}  
\usepackage{stmaryrd}  
\usepackage{amsmath,amsxtra,amsthm,amssymb,xr,enumerate,fullpage,comment, graphicx, mathtools}
\usepackage[all]{xy}
\usepackage[bbgreekl]{mathbbol}
\usepackage{bbm}
\usepackage{bm}
\usepackage{tikz-cd}
\makeatletter
\newlength{\doublefracgap}
\DeclareRobustCommand{\doublefrac}[2]{%
  \mathinner{\mathpalette\doublefrac@{{#1}{#2}}}%
}
\newcommand{\doublefrac@}[2]{\doublefrac@@#1#2}
\newcommand{\doublefrac@@}[3]{%
  \ooalign{%
    \raisebox{\doublefracgap}{$\m@th#1\frac{#2}{\phantom{#3}}$}\cr
    \raisebox{-\doublefracgap}{$\m@th#1\frac{\phantom{#2}}{#3}$}\cr
  }%
}

\makeatother

\DeclareMathSymbol{\invques}{\mathord}{operators}{`>}
\DeclareUnicodeCharacter{00BF}{\tmquestiondown}
\DeclareRobustCommand{\tmquestiondown}{%
  \ifmmode\invques\else\textquestiondown\fi
}

\usepackage{verbatim}
\numberwithin{equation}{section}

\makeatletter
\newcommand{\mylabel}[2]{#2\def\@currentlabel{#2}\label{#1}}
\makeatother

\newtheorem{theorem}{Theorem}[section]
\newtheorem{lemma}[theorem]{Lemma}

\newtheorem{proposition}[theorem]{Proposition}

\newtheorem{corollary}[theorem]{Corollary}
\newtheorem{defn}[theorem]{Definition}
\newtheorem{example}[theorem]{Example}

\newtheorem{remark}[theorem]{Remark}

\newtheorem{lthm}{Theorem}

\newtheorem{lcorollary}[lthm]{Corollary}
\newcommand{\rank}{\mathrm{rank}}

 \DeclareFontFamily{U}{wncy}{}
    \DeclareFontShape{U}{wncy}{m}{n}{<->wncyr10}{}
    \DeclareSymbolFont{mcy}{U}{wncy}{m}{n}
    \DeclareMathSymbol{\Sha}{\mathord}{mcy}{"58} 

\newcommand{\cS}{\mathcal{S}}

\newcommand{\length}{\mathrm{length}}
\newcommand{\ES}{\mathbf{ES}}

\newcommand{\pr}{\mathrm{pr}}
\newcommand{\cM}{\mathcal{M}}

\newcommand{\Tr}{\operatorname{Tr}}

\newcommand{\Gal}{\operatorname{Gal}}

\newcommand{\QQ}{\mathbb{Q}}
\newcommand{\Qp}{{\mathbb{Q}_p}}
\newcommand{\Zp}{\mathbb{Z}_p}
\newcommand{\ZZ}{\mathbb{Z}}
\renewcommand{\AA}{\mathbb{A}}
\newcommand{\FF}{\mathbb{F}}

\newcommand{\fF}{\mathfrak{F}}

\newcommand{\ord}{\mathrm{ord}}

\newcommand{\fp}{\mathfrak{p}}
\newcommand{\fq}{\mathfrak{q}}

\newcommand{\cH}{\mathcal{H}}
\newcommand{\cO}{\mathcal{O}}

\newcommand{\GL}{\mathrm{GL}}
\newcommand{\T}{\mathbb{T}}

\newcommand{\image}{\mathrm{Im}}
\newcommand{\cyc}{\textup{cyc}}

\newcommand{\Frob}{\mathrm{Frob}}

\newcommand{\fL}{\mathfrak{L}}

\newcommand{\Hom}{\mathrm{Hom}}
\newcommand{\Sel}{\mathrm{Sel}}
\newcommand{\Char}{\mathrm{char}}
\newcommand{\Ind}{\mathrm{Ind}}
\newcommand{\ac}{\textup{ac}}
\newcommand{\LL}{\Lambda}
\newcommand{\TT}{\mathbb{T}}

\newcommand{\f}{\textup{\bf f}}

\newcommand{\fn}{\mathfrak{n}}

\newcommand{\lra}{\longrightarrow}

\newcommand{\res}{\textup{res}}

\newcommand{\fP}{\mathfrak{P}}

\newcommand{\fM}{\m}

\newcommand{\Tw}{\mathrm{Tw}}

\newcommand{\re}{\mathrm{red}}

\newcommand{\red}[1]{\textcolor{red}{#1}}
\newcommand{\blue}[1]{\textcolor{blue}{#1}}

\newcommand{\q}{\mathfrak{q}}

\newcommand{\m}{\mathfrak{m}}

\newcommand{\cD}{\mathcal{D}}
\newcommand{\Frac}{\mathrm{Frac}}

\newcommand{\su}{{\mathscr U}}

\newcommand{\fN}{\mathfrak{N}}

\newcommand{\Fp}{\mathbb{F}_p}
\newcommand{\fa}{\mathfrak{a}}

\newcommand{\rec}{\mathrm{rec}}

\newcommand{\End}{{\rm End}}

\newcommand{\Q}{\QQ}
\newcommand{\Z}{\ZZ}

\newcommand{\C}{{\mathbb C}}
\newcommand{\R}{{\mathbb R}}

\usepackage{hyperref}

\definecolor{pinegreen}{rgb}{0.0, 0.47, 0.44}

 \definecolor{pAlgae}{RGB}{87,115,135}
\definecolor{airforceblue}{rgb}{0.36, 0.54, 0.66}
	\definecolor{bondiblue}{rgb}{0.0, 0.58, 0.71}
\definecolor{britishracinggreen}{rgb}{0.0, 0.26, 0.15}
\definecolor{camouflagegreen}{rgb}{0.47, 0.53, 0.42}
\definecolor{darkcyan}{rgb}{0.0, 0.55, 0.55}

\hypersetup{
    colorlinks = true,
    linkcolor=blue,
     filecolor=blue,
     citecolor = darkcyan,      
     urlcolor=cyan,
    linkbordercolor = {white},
}
\title{Anticyclotomic Iwasawa theory at potentially ordinary primes}
\author{Antonio Lei}
\address{Antonio Lei\newline Department of Mathematics and Statistics\\University of Ottawa\\
150 Louis-Pasteur Pvt\\
Ottawa, ON\\
Canada K1N 6N5}
\email{antonio.lei@uottawa.ca}

\author{Peter Neamti}
\address{Peter Neamti\newline UCD School of Mathematics and Statistics\\ University College Dublin\\ Belfield\\Dublin 4\\Ireland}
\email{peter.neamti@ucdconnect.ie}

\subjclass[2020]{Primary 11R23; Secondary 11F11, 11G05}
\keywords{Anticyclotomic Iwasawa theory, modular forms,
Heegner points, Selmer groups, potentially ordinary reduction.}

\begin{document}

\begin{abstract}
We study anticyclotomic Iwasawa theory for self-dual twists of $p$-ordinary modular forms of weight two with nebentype $\epsilon^2$, where $\epsilon$ is a non-trivial Dirichlet character of conductor $p$. Under a generalized Heegner hypothesis and suitable residual hypotheses, we prove one divisibility in the anticyclotomic main conjecture in both the definite and indefinite settings. The prime $p$ may split, remain inert, or ramify in the imaginary quadratic field. Our proofs use a compatible bipartite Euler system of Heegner classes and theta elements, together with a comparison with the $p$-adic $L$-functions of Kohen--Pacetti. In the definite case, we prove cotorsion of the Selmer group and the characteristic-ideal divisibility by the $p$-adic $L$-function. In the indefinite case, we prove a rank-one statement and the corresponding divisibility in terms of the Heegner class. We also discuss criteria for equality and their relation to Kolyvagin's conjecture. Specialization gives finite-level Bloch--Kato Selmer bounds with local terms, and applications to Mordell--Weil ranks and $p$-primary Tate--Shafarevich groups for elliptic curves that have potentially ordinary reduction at $p$ satisfying our hypotheses.
\end{abstract}

\maketitle
\section{Introduction}
Let $p>3$ be a fixed prime. This article proves one divisibility in the anticyclotomic Iwasawa main conjectures when $p$ is an additive potentially ordinary prime, treating both the definite and indefinite settings.
Bertolini--Darmon proved one divisibility in the definite case for elliptic curves at ordinary primes, including multiplicative reduction \cite{BertoliniDarmon2005}, and Howard developed the bipartite Euler-system formalism linking the definite and indefinite settings \cite{howard06}. Darmon--Iovita established the corresponding results for signed Selmer groups at supersingular primes split in the imaginary quadratic field \cite{darmoniovita}, while Chida--Hsieh extended the ordinary method to higher-weight modular forms \cite{chidahsiehanticyclomainconjformodformscomposito}.
More recently, signed anticyclotomic Iwasawa theories for both definite and indefinite settings, treating abelian varieties of $\mathrm{GL}_2$-type at split non-ordinary primes and elliptic curves at inert supersingular primes were studied in \cite{BBL1,BBL2}.
In the good ordinary setting, Burungale--Castella--Kim proved Perrin-Riou's Heegner point main conjecture under suitable hypotheses \cite{BCK}. Bertolini--Longo--Venerucci proved both definite and indefinite main conjectures for elliptic curves at good ordinary and supersingular primes unramified in the quadratic field under suitable hypotheses, with one divisibility in the exceptional plus-signed inert supersingular case \cite{BLV}. At additive primes, Kohen--Pacetti constructed anticyclotomic $p$-adic $L$-functions in \cite{kohen_pacetti_anticyc_l-function}. One of the motivations of the present work is to study the corresponding Iwasawa main conjecture in the additive setting. Recently, Keller--Yin proved anticyclotomic Iwasawa main conjectures for elliptic curves with potentially good ordinary reduction at Eisenstein primes \cite{Keller2024-jg}.
Their results concern the residually reducible case, whereas we assume residual absolute irreducibility in the present article.

Let $E/\Q$ be an elliptic curve of conductor $Np^2$, where $(N,p)=1$, that has potentially ordinary reduction at $p$. As explained in  \cite{delbourgo_Iwasawa_unstable, delbourgo_p-adic_BSD}, there exists a Dirichlet character $\epsilon$ of conductor $p$, such that the twist
$f_E \otimes \epsilon$
is ordinary at $p$, where $f_E$ is the weight two newform attached to $E$. This twisting method is utilized in \cite{delbourgo_Iwasawa_unstable, delbourgo_p-adic_BSD, muller_CM} to study the cyclotomic Iwasawa theory at potentially ordinary primes  (see also \cite{forras_lei} for an approach to the potentially supersingular case). In this article, we develop the corresponding theory over the anticyclotomic $\Z_p$-extension of an imaginary quadratic field $K$, allowing $p$ to split, remain inert or ramify in $K$. The main inputs are the bipartite Euler system machinery developed in \cite{BertoliniDarmon2005, howard06}, together with the anticyclotomic $p$-adic $L$-function constructed by Kohen-Pacetti \cite{kohen_pacetti_anticyc_l-function}. We impose the following (Heegner) hypothesis on $N$ for the rest of the article:

\begin{itemize}
    \item[(\mylabel{item_Heeg}{\textbf{Heeg}})] $N=N^+N^-$ with $N^-$ square-free. Every prime $q \mid N^+$ splits in $K$, and every prime $q \mid N^-$ is inert in $K$. 
\end{itemize}
\begin{comment}
We deduce one divisibility in the main conjecture for potentially ordinary elliptic curves as a corollary of a more general result which applies for any $p$-ordinary weight two newform 
$$f\in S_2(\Gamma_0(N) \cap \Gamma_1(p), \epsilon^2),$$ where $\epsilon^2 \neq \mathds 1$ is a Dirichlet character of conductor $p$ (in particular, the $p$-adic $G_{\Q_p}$-representation attached to such a form is not semistable. It is potentially crystalline, becoming crystalline as a $G_{\Q_p(\mu_p)}$-representation). Let $\widetilde{f}:= \Tw_{\epsilon}(f)$ denote the self-dual twist of $f$ given as in Remark \ref{remark_twisting_example}, $\Sel_{N^-}(K_\ac, A_{\widetilde{f}})$ and (respectively $\Sel_{N^-}(K_\ac, T_{\widetilde{f}})$) be the discrete (respectively compact) Selmer group defined in \S \ref{section_Local_cohomology_and_Selmer_groups}, $\mathscr{L}_f$ be the $p$-adic $L$-function constructed by Kohen-Pacetti, and let $\omega(M)$ denote the number of distinct prime factors of some integer $M$. Then we have the following result. 
\end{comment}

Our main results apply more generally to a normalized $p$-ordinary eigenform
\[
 f\in S_2(\Gamma_0(N)\cap\Gamma_1(p),\epsilon^2),
\]
that is new away from $p$. Here, $\epsilon$ is a non-trivial Dirichlet character of conductor $p$. If $\epsilon^2\ne\mathbb{1}$, then $f$ is a newform of level $Np$. Otherwise, $f$ is the ordinary $p$-stabilization of a newform $f^\circ\in S_2(\Gamma_0(N))$.

Let $\widetilde f=\operatorname{Tw}_\epsilon(f)$ be its self-dual twist. Write $\Sel_{N^-}(K_{\ac},A_{\widetilde f})$ and $\Sel_{N^-}(K_{\ac},T_{\widetilde f})$ for the discrete and compact Selmer groups defined in \S\ref{section_Local_cohomology_and_Selmer_groups}, respectively, and let $\mathscr L_f$ be the $p$-adic $L$-function of \S\ref{sec:p-adic-L}, obtained from the construction of Kohen--Pacetti. For a positive integer $M$, let $\omega(M)$ denote the number of its distinct prime divisors. 
The hypotheses \eqref{item_CR} and \eqref{item_irr} are stated in \S\ref{sec:admissible} and \S\ref{section_Local_cohomology_and_Selmer_groups}, respectively.
We prove the following results.

\begin{lthm}[Theorem \ref{thm_definite_main_conj}] \label{thmA}
    Assume the hypotheses \eqref{item_Heeg},  \eqref{item_CR} and  \eqref{item_irr} hold, and that $\omega(N^-)$ is odd. Then $\Sel_{N^-}(K_\ac, T_{\widetilde{f}}) = 0$, $\Sel_{N^-}(K_\ac, A_{\widetilde{f}})^\vee$ is $\Lambda(\Gamma_\ac)$-torsion, and we have the divisibility 
    $$\left.\Char_{\Lambda(\Gamma_\ac)}\left(\Sel_{N^-}(K_\ac, A_{\widetilde{f}})^\vee \right) \right| (\mathscr{L}_f).$$

\end{lthm}

\begin{lthm}[Theorem \ref{them_indef_main_conjecture}]\label{thmB}
    Assume the hypotheses \eqref{item_Heeg},  \eqref{item_CR} and  \eqref{item_irr} hold, and that $\omega(N^-)$ is even. Then $\Sel_{N^-}(K_\ac, T_{\widetilde{f}})$ is a torsion-free $\Lambda(\Gamma_\ac)$-module and 
    $${\rm rank}_{\Lambda(\Gamma_\ac)} \left(\Sel_{N^-}(K_\ac, T_{\widetilde{f}})\right) = {\rm rank}_{\Lambda(\Gamma_\ac)} \left(\Sel_{N^-}(K_\ac, A_{\widetilde{f}})^\vee\right) =1.$$ Moreover, we have the divisibility
    $$\left.\Char_{\Lambda(\Gamma_\ac)}\left(\Sel_{N^-}(K_\ac, A_{\widetilde{f}})_{\rm tor}^\vee \right)\right|\Char_{\Lambda(\Gamma_\ac)}\left(\Sel_{N^-}(K_\ac, T_{\widetilde{f}})/ \kappa_{\infty} \right)^2 ,$$
    where $\kappa_\infty$ is the Heegner class defined in \eqref{eq:kappa_infty}.
\end{lthm}
\begin{comment}
\blue{Maybe we can remove this remark now?}
\red{\begin{remark}
    If $\epsilon^2 = \mathbb 1$ and $a_p(f)^2 \not \equiv 1 \mod \varpi$, where $\varpi$ is a uniformizer for the $p$-adic Hecke field of $f$, then $f$ is the $p$-stabilization of a newform $f^\circ \in S_2(\Gamma_0(N))$. If $\epsilon = \mathbb 1$ and $(D_K, p) =1$, then this case is handled by \cite{chidahsiehanticyclomainconjformodformscomposito}. Building on their method, the arguments in this paper extend to the case when $\epsilon$ is the quadratic character of conductor $p$ and also allow for $p \mid D_K$. The condition that $a_p(f)^2 \not \equiv 1 \mod \varpi $ ensures that Lemma \ref{lemma_vanishing_cohomology_at_p} still holds, and (more importantly) it is necessary to make use of the version of Ihara's lemma given in \cite[\S 5.2]{chidahsiehanticyclomainconjformodformscomposito}.
\end{remark}}
\end{comment}

%Taking $f = f_E \otimes \epsilon$ under the assumption that $\epsilon^2\ne\mathbb{1}$, we see that the $p$-adic Galois representation attached to $\widetilde{f}$ is isomorphic to the $p$-adic Galois representation attached to $E$ (\textit{cf.} Remark \ref{remark_twisting_example}). We then obtain the following as an immediate corollary of the results above.

For an elliptic curve $E/\QQ$ as above, take $f$ to be the $p$-ordinary eigenform associated with $f_E\otimes\epsilon$ as specified above.
After extending coefficients, the representation attached to $\widetilde f$ identifies with that on the Tate module $T_p(E)$.
In this setting, \eqref{item_CR_2} is the absolute irreducibility of $E[p]$; we discuss the interpretation of the remaining hypotheses in Remark~\ref{rk:EC_conditions}.
Theorems~\ref{thmA} and \ref{thmB} give the following corollary.

\begin{lcorollary}
    Assume that the hypotheses \eqref{item_Heeg}, \eqref{item_CR} and  \eqref{item_irr} hold for the $p$-ordinary eigenform $f$ specified above. 
    \begin{itemize}
        \item[i)] If $\omega(N^-)$ is odd, then $\Sel(K_\ac, E[p^\infty])^\vee$ is $\Lambda(\Gamma_\ac)$-torsion and the divisibility of Theorem A holds with $A_{\widetilde{f}}$ replaced by $E[p^\infty]$.
        \item[ii)] If $\omega(N^-)$ is even, then ${\rm rank}_{\Lambda(\Gamma_\ac)} \left(\Sel(K_\ac, E[p^\infty])^\vee\right) =1$ and the divisibility of Theorem B holds with $A_{\widetilde{f}}$ replaced by $E[p^\infty]$.
    \end{itemize}
\end{lcorollary}

To prove Theorems~\ref{thmA} and \ref{thmB}, we construct  a bipartite Euler system for the self-dual twist $\widetilde{f}$, relating Heegner classes to the theta elements of \cite{kohen_pacetti_anticyc_l-function}.
Although we use the twisting method of \cite{delbourgo_Iwasawa_unstable,delbourgo_p-adic_BSD}, the passage to anticyclotomic Iwasawa main conjectures is not formal.
When $\epsilon^2\ne\mathbb{1}$, the norm pullback of $\epsilon$ to $K$ is not anticyclotomic, so the twist cannot be absorbed into the anticyclotomic character variable.
We therefore incorporate the twist into the constructions of Heegner classes and theta elements and establish their compatibility under reduction modulo powers of $\varpi$.
Hida theory supplies the required integral level control and, in the case where $\epsilon^2\ne\mathbb{1}$, a comparison modulo $\varpi$ with higher-weight specializations.
In the indefinite case, we also prove the non-vanishing modulo $\varpi$ of the anticyclotomic Heegner class, which is the key ingredient in the proof of the rank-one statements in Theorem~\ref{thmB}.

In \S\ref{s:kolyvagin}, we discuss the relation between Kolyvagin's conjecture and the criteria for equality in Theorems~\ref{thmA} and \ref{thmB}. In \S\ref{s:app}, we discuss arithmetic applications of these theorems. 
In the definite case, we show that $L(\widetilde f/K,1)\ne0$ implies the finiteness of the Bloch--Kato Selmer group, whose length, together with local Tamagawa terms, is bounded by the valuation of $\mathscr L_f$ at the trivial character.
In the indefinite case, a nonzero specialized Heegner class gives a Selmer group of corank one and bounds the torsion in its Pontryagin dual, together with a finite local-condition correction, in terms of the index of the class.

\subsubsection*{Acknowledgements} 
The second-named author (P.N.) wishes to thank the first-named author (A.L.) for agreeing to host him at the University of Ottawa during the Summer of 2026 as part of the University of Ottawa's `Doctoral Mobility Program'. This work was first started during this visit. P.N. is also supported by a NUI travelling studentship. The research of A.L. is supported by the NSERC Discovery Grants Program RGPIN-2026-07384.
A.L. also thanks Daniel Delbourgo and Ariel Pacetti for interesting discussions on topics related to this article.

\section{Preliminaries}
Let $N \in \Z$, $p>3$ be a fixed prime, and $f\in S_2(\Gamma_0(N) \cap \Gamma_1(p), \epsilon^2) $ all be as in the introduction.

Fix embeddings 
$$i_{\C}: \overline \Q \hookrightarrow \C, \quad i_{\overline{\Q}_p}: \overline \Q \hookrightarrow \overline{\Q}_p, $$ 
and let $v_p$ be the induced $p$-adic valuation, normalized such that $v_p(p)=1$. We impose the following ordinarity condition throughout the article.
\begin{itemize}
    \item[(\mylabel{item_ord}{\textbf{ord}})] $v_p(a_p(f)) = 0$. 
\end{itemize}
If $\chi: (\Z/M\Z)^\times \lra \C^\times$ is a Dirichlet character, then the projection $\AA^\times \simeq \Q^\times \times \widehat \Z^\times \times  \R_{> 0} \twoheadrightarrow (\Z/M\Z)^\times$ allows us to pull $\chi$ back to an idelic Hecke character, $\chi_{\AA^\times}$, satisfying 
\begin{equation}
    \chi_{\AA^\times}(q\cdot \alpha \cdot r) = \chi^{-1}(\overline \alpha), \quad q\in \Q^\times, \ \alpha \in \widehat \Z^\times, \ r \in \R_{> 0},\label{eqn_Dirichlet_to_Hecke}
\end{equation}
where $\overline \alpha$ is the image of $\alpha$ under the projection $\widehat \Z^\times \twoheadrightarrow (\Z/M\Z)^\times$.

\subsection{Iwasawa Algebras}
Let $K$ be an imaginary quadratic field, and write $D_K=|\operatorname{disc}(K)|$.
 Let $H_\infty:= \bigcup_n H_n$, where $H_n$ is the ring class field of $K$ of conductor $p^n$. Let $\Delta_\ac$ denote the torsion subgroup of $\Gal(H_\infty/K)$, so that we have the split short exact sequence 
$$0\lra \Delta_\ac \lra \Gal(H_\infty/K) \lra \Gamma_\ac \lra 0,$$
where $\Gamma_\ac \cong \Zp$ as pro-$p$ abelian groups. Let $K_\ac$ denote the anticyclotomic $\Zp$-extension of $K$, so that $\Gal(K_\ac/ K) = \Gamma_\ac$, and let $K_m\subset K_\ac$ be the field extension of $K$ given by {$K_m := H_m \cap K_\ac$}. 
If $\mathcal{O}$ is some $p$-adically complete $\Zp$-algebra (most commonly for us, the ring of integers of some finite extension of $\Qp$), then we define the Iwasawa algebra
$$\LL_\mathcal{O}(\Gamma_\ac) := \mathcal{O} \llbracket \Gamma_\ac \rrbracket.$$
If $T$ denotes a transcendental variable, then a choice of topological generator $\gamma_\ac \in \Gamma_\ac$ induces an isomorphism of topological rings 
$$\Lambda_\mathcal{O}(\Gamma_\ac) \xrightarrow{\gamma_\ac \mapsto 1+T} \mathcal{O}[[T]]. $$
Unless there is danger of confusion, we will often drop the $\mathcal{O}$ from the notation and just write $\Lambda(\Gamma_\ac)$.

\section{Quaternion algebras and level-raising}
\subsection{Quaternion algebras and modular forms}
\subsubsection{Quaternion algebra preliminaries}\label{subsub_QuaternionPrelims}
Let $B$ denote a quaternion algebra over $\Q$ which is either ramified at $\infty$ (definite), or split at $\infty$ (indefinite). Let $\Delta \in \Z$ denote the discriminant of $B$, and assume it is coprime to $p$. We fix an Eichler order $R_{M^+} \subset B$ of level $M^+$, and write $\widehat B := B \otimes \widehat \Z$, $\widehat R_{M^+} := R_{M^+} \otimes \widehat \Z$ for the respective profinite completions, and $B_v := B \otimes \Q_v$, $R_{M^+, v} : = R_{M^+} \otimes \Z_v$. If $b \in \widehat B$, then we write $b_v$ for the image of $b$ in $B_v$ under the natural map. By definition of $B$, there exists an isomorphism
$$\iota: \prod_{v \nmid \Delta}B_v \cong \prod_{v\nmid \Delta}M_2(\Q_v), $$
which we henceforth fix so that each maximal order $R_{M^+, v} \subset B_v$, $v\nmid M^+ \Delta$, is identified with $M_2(\Z_v)$. For $v \nmid \Delta$, let $\iota_v: B_v \cong M_2(\Q_v)$ be the induced isomorphism on local components. Note that $\iota$ induces an isomorphism 
$$\iota: \prod_{v\nmid \Delta}R_{M^+, v} \cong \prod_{v\nmid \Delta M^+}M_2(\Z_v) \times \prod_{v \mid M^+}\left \{\begin{pmatrix}
    a &b \\ c & d 
\end{pmatrix} \in M_2(\Z_v) :c \equiv 0 \mod v^{\ord_v(M^+)} \right\}.$$
For $n\geq 1$, let $\su_n \subset \widehat R_{M^+}^\times$ denote the open compact subgroup of $\widehat B^\times$ given by 
\begin{align}
\mathscr{U}_n := \{r \in \widehat R_{M^+}^\times : \iota_p(r_p) \equiv \begin{pmatrix}
    a & b \\
    0 & a 
\end{pmatrix} \mod p^n\}. \notag
\end{align}

\begin{remark}\label{remark_Comparing_level_subgroups}
    Let $\mathscr{U}^1_n := \{r \in \widehat R_{M^+}^\times : \iota_p(r_p) \equiv \begin{pmatrix}
    a & b \\
    0 & 1 
\end{pmatrix} \mod p^n\}.$
    By construction, $ \su^1_n, \mathscr U_n \subset \mathscr U^0_n : = \{r \in \widehat R_{M^+}^\times: \iota_p(r_p) \equiv \begin{pmatrix}
        a & b \\
        0 & d
    \end{pmatrix} \mod p^n \}$, and we have isomorphisms
    \begin{align}
        \mathscr U^0_n / &\mathscr U^1_n \cong (\Z/ p^n\Z)^\times \notag, 
        &\begin{pmatrix}
            a & b \\
            c & d
        \end{pmatrix} \mod \su^1_n\mapsto {d}, \notag \\
        \mathscr U^0_n / &\mathscr U_n \cong (\Z/ p^n\Z)^\times \notag, 
        &\begin{pmatrix}
            a & b \\
            c & d
        \end{pmatrix} \mod \su_n \mapsto \frac{a}{d}. \notag 
    \end{align}  
    Let $A$ be a ring containing $\mu_{p^{n-1}(p-1)}$, and such that if $n\neq1$, then $A[\frac{1}{p}] = A$. Let $ \mathscr{M}$ be an $A[\mathscr U^0_n]$-module, $\chi_1, \chi_2$ be two characters of $(\Z/p^n\Z)^\times$, and let us denote by $ \mathscr{M}(\chi_1, \chi_2)$ the subspace of $ \mathscr{M}$ on which all $ \begin{pmatrix}
        a & b\\
        c & d
    \end{pmatrix}\in \mathscr U^0_n$ act by $\chi_{1, \AA^\times}(a)\cdot \chi_{2, \AA^\times}(d)$ (\textit{cf.} \eqref{eqn_Dirichlet_to_Hecke} for idelic character convention). In particular, we have the decompositions 
        \begin{align}
            & \mathscr{M}^{\mathscr U^1_n} = \bigoplus_\chi \mathscr{ \mathscr{M}}(\mathbb 1, \chi), \label{eqn_FirstChoice_decomp} \\
            & \mathscr{M}^{\mathscr U_n}  = \bigoplus_{\chi}  \mathscr{M}(\chi^{-1}, \chi),\label{eqn_SecondChoice_decomp}
        \end{align}
    where the sum is over all characters of $(\Z/p^n\Z)^\times$. It is clear that the central character at $p$ of any element in $ \mathscr{M}^{\mathscr{U}_n}$ is trivial, making this choice of subgroup natural for considering forms which have self-dual Galois representations.
    Suppose now that for a Dirichlet character $\chi$, we have an isomorphism
    $$\Tw_{\chi}: \mathscr M \xrightarrow{ \sim} \mathscr M \otimes \chi_{\AA^\times} \circ \det.$$
    Comparing \eqref{eqn_FirstChoice_decomp} with \eqref{eqn_SecondChoice_decomp}, we see that the map $\Tw_{\chi}$ restricts to an isomorphism of $A$-modules
    $$\Tw_{\chi} : \mathscr{M}(\mathbb{1}, \chi^{2})
    \xrightarrow{\ \sim\ } \mathscr{M}(\chi^{-1}, \chi),$$
    which satisfies
    $\Tw_{\chi}(u \cdot f) = \chi_{\mathbb{A}^\times}(\det u)\,
    \bigl(u \cdot \Tw_{\chi}(f)\bigr)$ for $f\in \mathscr M(\mathbb 1, \chi^2)$, $u \in \su_n^0$.
\end{remark}

\begin{remark} \label{remark_twisting_example}
    Suppose momentarily that $\Delta=1$ and $M^+=N$. Let $M^1_n$ denote the space of weight $2$ modular forms of level $\su^1_n$ and let $M_n$ denote the weight 2 forms of level $\su_n$ (\textit{cf.} Definition \ref{def_indef-analytic}). Let $E$ be an elliptic curve of conductor $Np^2$ which has potentially ordinary reduction at $p$ as in the introduction. Let $f_E \in M^1_2(\mathbb 1,\mathbb 1)$ denote the associated weight $2$ form. Let $\epsilon$ be the Dirichlet character of conductor $p$ as in the introduction, and let $f\in M^1_1(\mathbb 1, \epsilon^2)$ be the $p$-ordinary eigenform as specified in the introduction. If $M_1$ denotes the weight $2$ forms invariant under $\su_1$, then there is another twisting map given by the isomorphism of $\C$-vector spaces 
    \begin{align*}
        \Tw_{\epsilon}: M^1_1(1, \epsilon^{2}) &\xrightarrow{\sim} M_1(\epsilon^{-1}, \epsilon) \\
        \Tw_{\epsilon}(g)(z, \gamma) = \epsilon_\infty({\rm Im}(z))\cdot\epsilon^{-1}_f ( \det\gamma)\cdot &g(z, \gamma), \quad (z, \gamma) \in \mathcal{H}^\pm \times \GL_2(\AA_f),
    \end{align*}
    where we decompose $\epsilon = \epsilon_\infty \cdot \epsilon_f$ into its infinite and finite parts respectively (note that $\epsilon_\infty $ depends only on the sign of its input).
    Then $\Tw_{\epsilon}(f)=: \widetilde f \in M_1(\epsilon^{-1}, \epsilon)$, whereas $f_E\in M^1_2(\mathbb 1, \mathbb 1)$, so the two twisting maps are not inverses of each other; however, the Galois representations attached to $f_E$ and $\widetilde f$ are isomorphic.
\end{remark} 

\subsubsection{Shimura varieties for quaternion algebras}\label{s:components}
Suppose that $B$ is a definite or indefinite quaternion algebra over $\Q$, and $U \subset \widehat B^\times$ is a compact open subgroup. Set $\mathcal{H}^\pm = \C \setminus \R$, and define the (analytification of the) Shimura set/curve in the definite/indefinite case by
$$X(U)(\C) := \begin{cases}
    B^\times \backslash \widehat B^\times /U, &\text{ if } B \text{ is definite,} \\
    B^\times \backslash \mathcal{H}^\pm \times \widehat B^\times /U, &\text{ if } B \text{ is indefinite.}
\end{cases}$$
If the discriminant equals $1$, then we instead consider the compactification of $X(U)(\C)$ obtained by adjoining finitely many cusps. In the definite case, the left action of $B^\times$ on $\widehat B^\times$ is given by left multiplication through the embedding $B^\times \hookrightarrow \widehat B^\times$. In the indefinite case, $B^\times$ acts diagonally on $\widehat B^\times \times \mathcal{H}^\pm$, where the action on $\widehat B^\times$ is as before, and the action on $\mathcal{H}^\pm$ is through the embedding $\iota_\infty: B^\times \hookrightarrow \GL_2(\R)$, and the standard Moebius transformation. On the other hand, the right action of $U$ is only on $\widehat B^\times$, and is just given by right multiplication in both the definite and indefinite case. The Shimura set/curve $X(U)(\C)$ admits a canonical model over $\Q$, which we denote by $X(U)$. If $(z,x) \in \mathcal{H}^\pm \times\widehat B^\times$, then we denote by $[z, x]_U$ (resp. $[x]_U$) the image of $(z, x)$ (resp. $x$) in $X(U)(\C)$ if $B$ is indefinite (resp. definite).

In the definite case, $X(U)$ is a zero-dimensional variety comprised of finitely many points. In the indefinite case, $X(U)$ is a smooth projective curve of dimension one satisfying
$$\pi_0(X(U)_{\overline \Q}) \cong \Q^{\times}_{>0}\backslash\mathbb{A}_{f}^{\times}/\mathrm{N}(U),$$
where $\mathrm N:\widehat B^\times\longrightarrow\mathbb A_f^\times$ denotes the reduced norm.
In particular, $\mathrm{N}(\su_{n,p}) = (\Z_p^\times)^2$, so $X(\su_n)$ has \emph{two} connected components.

\subsubsection{Weight 2 forms for definite quaternion algebras}
For this subsubsection, $B$ is a \emph{definite} quaternion algebra over $\Q$ of discriminant $\Delta$. In particular, $\omega(\Delta)$ is \emph{odd}. The main reference for this subsection is \cite[\S 2.3, 4.1]{chidahseihanticycloLvaluescrelle}, although much of the discussion in \textit{op. cit.} is simplified greatly for weight $2$ forms.
\begin{defn}\label{def_definite_modularforms}
Let $A$ be a ring, and $U \subset \widehat B^\times$ an open compact subgroup. We denote by $M_2(U, A)$ the space of modular forms of level $U$ defined over $A$, consisting of functions $f : \widehat B^\times \lra A$ such that
$$f(\alpha b u) = f(b) \quad \text{for all } \alpha \in B^\times,\, u \in U.$$
\end{defn}
By definition, $$M_2(U, A) \simeq \Hom(X(U), A) \simeq H^0(X(U), A). $$
If $\psi$ is some Dirichlet character of conductor ${\rm cond}(\psi)$, we may view it as a Hecke character or a $p$-adic Hecke character via our fixed embeddings $i_\C$, $i_{\overline{\Q}_p}$, and, abusing notation, we will continue to denote both incarnations by $\psi$. If $C$ is a subring of $\C$ or $\overline \Q_p$ such that $C = C[\frac{1}{\rm{cond}(\psi)}]$, then we set 
$$M_2(U, C, \psi):= \{f \in M_2(U, C): f(\alpha b) = \psi(\alpha)f(b), \quad \alpha \in \AA_f^\times \subset \widehat B^\times\}.$$

\subsubsection{Weight 2 forms for indefinite quaternion algebras}
For this subsubsection, $B$ is an \emph{indefinite} quaternion algebra over $\Q$, of discriminant $\Delta$, so that $\omega(\Delta)$ is even.

First, if $\gamma \in B^\times \hookrightarrow B^\times \otimes \R$, is such that $\iota_\infty(\gamma \otimes 1) =  \begin{pmatrix}
    a & b\\
    c & d
\end{pmatrix},$ then set 
$$j(\gamma,z) = cz+d.$$
\begin{defn}\label{def_indef-analytic}
    Let $U \subset \widehat B^\times$ be an open compact subgroup. Set $M_{2}(U,\C)$ to denote the space of functions
    $f: \mathcal{H}^{\pm}\times\widehat{B}^{\times}\longrightarrow\C$ such that
    \begin{enumerate}
      \item[i)] $z\mapsto f(z,b)$ is holomorphic for every $b\in\widehat{B}^{\times}$;
      \item[ii)] $f(\gamma z,\gamma b u)=\det(\gamma)^{-1}j(\gamma,z)^{2}\,f(z,b)$ for all $\gamma\in B^{\times}$ and $u\in U$.
      \item[iii)] If $\Delta=1$, then for all $b\in \widehat B^\times$, $f(\cdot, b)$ is holomorphic at the cusps of $\Gamma_b: = bUb^{-1} \bigcap \GL_2(\Q)^+$.
    \end{enumerate}
\end{defn}

\subsubsection{Hecke operators}\label{S:hecke}
For this subsubsection, $B$ is either definite or indefinite, and $U \subset \widehat B^\times $ is an open compact subgroup. Let $\mathcal{O}$ denote the ring of integers of a finite extension of $\Q_p$. If $x \in \widehat B^\times$, then we have the double coset decomposition 
\begin{equation}
U x U = \bigsqcup_{i\in I} x_iU,\label{eqn_double_coset_decomp}
\end{equation}
for some finite index set $I$. Let $[U x U] \in \Hom_{\mathcal{O}}(M_2(U, \mathcal{O}), M_2(U, \mathcal{O}))$ be the Hecke operator defined by 
$$[U x U] f(z,b) = \sum_{i}f(z,bx_i), $$
where in the definite case, we are abusing notation in setting $f(z, b): =f(b)$.
If $v$ is a finite place of $\Q$ and $v\nmid \Delta$, then we fix elements 
\begin{align} 
    &\pi_v = \begin{pmatrix}
        v & 0\\
        0 & 1
    \end{pmatrix},  \ z_v = \begin{pmatrix}
        v & 0\\
        0 &v
    \end{pmatrix} \in \GL_2(\Q_v).  \notag
\end{align}
If $v\mid \Delta$, then we fix an element 
$$ \pi'_v \in B_v^\times: \ N(\pi'_v)=v,  \notag $$
where $N$ is the reduced norm of $B$. Then we can define the Hecke algebra $\T_{\Delta}$ as the subalgebra of $\End_{\mathcal O}(M_2(\su_n, \mathcal O))$ generated by the $\{T_v,S_v, S_v^{-1}: \ v\nmid \Delta M^+p\}\bigcup \{U_v: \ v\mid \Delta M^+p\} \bigcup \{\langle a \rangle: \ a \in \Z_p^\times\}$, where the operators $T_v,S_v, U_v$ are given by
\begin{align}
    &T_v := [\su_n \iota^{-1} (\pi_v) \su_n], \notag \\
    &S_v := [\su_n \iota^{-1} (z_v) \su_n],  \notag \\
    &U_v: =\begin{cases}
        [\su_n \iota^{-1} (\pi_v) \su_n], \ v \mid M^+p, \\
        [\su_n \pi'_v \su_n], \ v \mid \Delta,
    \end{cases} \notag
\end{align}
and the diamond operators for $a \in \Z_p^\times$ at $p$ are given by 
\begin{align}
    \langle a \rangle := [\su_n \iota^{-1}(\mathbf d(a))\su_n], && \mathbf d(a) = \begin{pmatrix}
        a &0 \\
        0 & 1
    \end{pmatrix} \in \GL_2(\Z_p). \notag
\end{align}
\begin{remark}\label{remark_S_q_trivial}
Note that the $S_v$ act trivially on $M_2(\su_n, \mathcal{O})$, but we include them in the definition for completeness.
\end{remark}

\begin{remark}\label{rk:Hecke}
    Recall the open compact subgroup $\su^1_n$ of Remark \ref{remark_Comparing_level_subgroups}. We can define the Hecke algebra for $\mathscr{U}^1_n$ in the same way as for $\su_n$. Suppose that $g\in M_2(\su^1_n, \C, \chi^{2})$ is a Hecke eigenform with eigensystem $(\lambda_q(g))_q$ and recall the twisting operator $\Tw_{\chi}$ of Remark \ref{remark_twisting_example} (which we may extend to the definite case by dropping the $\epsilon_\infty(\rm{Im}(z))$ term). Then $\Tw_{\chi}$ maps $g$ to the Hecke eigenform $\Tw_{\chi}(g)=: \widetilde g \in M_2(\su_n,\C)$, whose Hecke eigensystem $(\lambda_q(\widetilde g))_q$ is given by 
    \begin{align*}
        \begin{aligned}
        \lambda_q(\widetilde g)& = \chi^{-1}(q)\lambda_q(g), \quad q\nmid p \\
        \lambda_p(\tilde g) &= \lambda_p(g),  \\
        \langle a \rangle& = \chi(a), \quad a\in\Z_p^\times 
        \end{aligned}
    \end{align*}
\end{remark}

\subsection{Admissible forms}
\subsubsection{}\label{sec:admissible}
Recall the fixed eigenform $f \in S_2(\Gamma_0(N)\cap \Gamma_1(p), \epsilon^2)$, and the decomposition $N = N^+N^-$ (\textit{cf.} \eqref{item_Heeg}), where $N^-$ is a square-free product of primes inert in the fixed imaginary quadratic field $K$, and $(N^-, pN^+)=1$. From now on $B$ will always denote the fixed quaternion algebra of discriminant $N^-$, and $ \su_n \subset \widehat R_{N^+}^\times \subset \widehat B^\times$ are defined as in \S \ref{subsub_QuaternionPrelims}, with $M^+=N^+$. Let $\mathcal{O}$ denote the ring of integers of $\Q_p((a_n(f))_n)$, and $\varpi$ be a uniformizer of $\mathcal{ O}$. Since $\epsilon$ has conductor $p$, we can write $\epsilon = \omega^{k}$ and hence $\epsilon^2 = \omega^{2k}$ for some integer $k$ such that $0\le k\le p-2$, where $\omega$ denotes the Teichm\"uller character. We impose the following hypotheses for the rest of this article:
\begin{itemize}
    \item[(\mylabel{item_CR}{\textbf{CR}})] The following conditions hold:
    \begin{itemize}
        \item[(\mylabel{item_CR_1}{\textbf{CR 1}})] One of the following holds:
        \begin{itemize}
            \item[i)] $2k < p-3$ and $\# (\FF_p^\times)^{2k+1}>5$.
           \item[ii)] $k = \frac{p-1}{2}$, $p\geq 7$, and $a_p(f)^2 \not \equiv 1 \mod \varpi $.
        \end{itemize} 
        \item[(\mylabel{item_CR_2}{\textbf{CR 2}})] $\overline \rho_f$ is absolutely irreducible.
        \item[(\mylabel{item_CR_3}{\textbf{CR 3}})] If $\ell$ is a prime, then $\overline \rho_f$ is ramified at $\ell$ whenever 
        \begin{itemize}
            \item[(i)] $\ell \mid N^-$  and $\ell^2 \equiv 1 \mod p$,
            \item[(ii)] $\ell \mid N^+$.     
        \end{itemize}
        \item[(\mylabel{item_CR_4}{\textbf{CR 4}})] {$\overline\rho_f|_{I_\ell}$ is irreducible if $\ell^2|N^+$ and $p|\ell-1$.}
         \end{itemize}
\end{itemize}

The condition \eqref{item_CR_4} originates from \cite{CH-erratum}.

\begin{defn}[$n$-admissible primes]\label{def_n-admissible_primes}
    We say that a prime $\ell$ is $n$-admissible if the following conditions hold
    \begin{itemize}
    \item[i)] $\ell \nmid pN$.
    \item[ii)] $\ell$ is inert in the imaginary quadratic field $K$.
    \item [iii)] $p \nmid \ell^2-1$.
    \item[iv)] $(\ell+1)^2 \equiv \epsilon^{-2}(\ell)a_\ell(f)^2 \mod \varpi^n$.
    \end{itemize}
    
\end{defn}

Let $\Delta \in \Z$ be such that $(\Delta, N^+)=1$, and denote by $B_\Delta$ a quaternion algebra over $\Q$ of discriminant $\Delta$. If $B_\Delta$ is definite, we write $M_2^\Delta(\su_n, \mathcal{O}_n)$ for the space of $p$-adic modular forms of level $\su_n$ with coefficients in $\mathcal{O}_n:=\mathcal{O}/\varpi^n$ (\textit{cf.} Definition \ref{def_definite_modularforms}). If $B_\Delta$ is indefinite, then we set $M_2^\Delta(\su_n, \mathcal{O}_n):= \Hom(\T_{\Delta}, \mathcal{O}_n)$, where here $\Hom$ means ring homomorphism. Moreover, if $g \in M_2^\Delta(\su_n, \mathcal{O}_n)$, then we set $\lambda_q(g):= g(T_q)$ if $q\nmid \Delta N^+p$, and $\lambda_q(g):= g(U_q)$ if $q\mid \Delta N^+p$. We call elements of $M_2^\Delta(\su_n, \mathcal{O}_n)$ $\T_{\Delta}$-eigenforms.

\begin{defn}[$n$-admissible forms]\label{def_n-admissible_forms} 
    An $n$-admissible form is a tuple $(\Delta, g)$, where $\Delta \in \Z$ and $g \in M^\Delta_2(\su_n, \mathcal{O}_n)$ is a $\T_{\Delta}$-eigenform such that
    \begin{itemize}
    \item[i)] $N^- \mid \Delta$ and $\frac{\Delta}{N^-}$ is a product of $n$-admissible primes.
    \item[ii)] $g \not\equiv 0 \mod \varpi$.
    \item[iii)] We have the following congruences
    \begin{align*}
        &\lambda_q(g)\equiv \epsilon^{-1}(q)a_q(f) \mod \varpi^n, \ (q\nmid p\Delta), && \lambda_p(g)\equiv a_p(f) \mod \varpi^n,\\
        &\lambda_\ell(g) \equiv \varepsilon_\ell \mod \varpi^n, \ (\ell \mid \frac{\Delta}{N^-}), &&\langle a \rangle \equiv \epsilon(a) \mod \varpi^n, 
    \end{align*}
    where $\varepsilon_\ell = \pm 1$ is uniquely determined by the condition $\epsilon^{-1}(\ell) a_\ell(f)\equiv\varepsilon_\ell\cdot(\ell+1)\mod\varpi^n$.
    \end{itemize}
    
\end{defn}

\begin{example}\label{example_n-admissible_form}
    Let $f_B \in M_2(\su^1_1, \mathcal{O}, \epsilon^2) \subset M_2(\su^1_n, \mathcal{O}, \epsilon^2)$ be the modular form for the definite or indefinite quaternion algebra $B$ associated to $f$ by Jacquet--Langlands. Then $\Tw_\epsilon(f_B)\otimes \mathcal{O}_n =: \widetilde f_n \in M_2(\su_n, \mathcal{O}_n)$ and $(N^-, \widetilde f_n)$ is an $n$-admissible form.

    If $\epsilon^2=\mathbb{1}$, let $f_B^\circ$ be the Jacquet--Langlands transfer of $f^\circ$, where $f^\circ$ is the newform of level $N$ whose ordinary $p$-stabilization is $f$. Then $f_B$ is the ordinary stabilization of $f_B^\circ$.
\end{example}

\begin{remark} \label{remark_Hida_theory}
    Let $k_0 = 2k+2$ and assume that \eqref{item_CR_1}(i) holds. This is the condition $({\rm CR}^+)(1)$ in \cite{chidahsiehanticyclomainconjformodformscomposito} for a weight $k_0$ eigenform. If we consider the ideals 
    \begin{align*}
         &P_{\epsilon}=(\langle a\rangle-\epsilon(a):a\in\Zp^\times)\subset \TT_{\Delta}, & P_{k_0}=(\langle a\rangle-a^k:a\in\Zp^\times)\subset \TT_{\Delta},
    \end{align*}
    then 
    $$P_{\epsilon} \equiv P_{k_0} \mod \varpi.$$ 
    Moreover, if $\mathfrak{m}$ denotes the maximal ideal of $\TT_\Delta$ attached to an $n$-admissible form $\mathcal{D}$, then we can make use of Hida theory as in \textit{op. cit.} to relate the localization at $\mathfrak{m}$ of the ordinary subspace of weight two forms of level $\su_n$ modulo $\varpi$ annihilated by $P_\epsilon$ to the corresponding localized space of ordinary weight $k_0$ forms with level prime to $p$.

    If \eqref{item_CR_1}(ii) holds, then $\epsilon^2=\mathbb{1}$. Undoing the twist of Remark~\ref{rk:Hecke} reduces the comparison to weight two and trivial nebentype. Since $a_p(f)^2\not\equiv1\mod \varpi$, the operator $U_p^2-1$ is a unit at the corresponding maximal ideal. As it annihilates the $p$-new contribution, that contribution vanishes after localization. Thus, the arguments in \cite[Proposition~4.2, Lemma~5.3 and Corollary~5.4]{chidahsiehanticyclomainconjformodformscomposito} apply. Twisting back gives the required statements for the $\epsilon$-component that is being studied in the present paper.
\end{remark} 

\section{{Local cohomology and Selmer groups}}\label{section_Local_cohomology_and_Selmer_groups}

Let $\Delta \in \Z$ be a square-free product of primes such that $(\Delta, N^+)=1$, $N^- \mid \Delta$, and $\frac{\Delta}{N^-}$ is a product of $n$-admissible primes (\textit{cf.} Definition \ref{def_n-admissible_primes}). Fix $B_\Delta$ to be the quaternion algebra over $\Q$ of discriminant $\Delta$, and an Eichler order $R=R_{N^+} \subset B_\Delta$ of level $N^+$ as in \S \ref{subsub_QuaternionPrelims}. We write $\mathcal{O}_K[p^m] := \Z + p^m \mathcal{O}_K$ for the order of $K$ of conductor $p^m$, and denote the ring class field of conductor $p^m$ by $H_m$. We also impose the following hypothesis: 
\begin{itemize}
    \item[(\mylabel{item_irr}{\textbf{irr}})]  If $N^- =1$, then $p\mid \#\overline \rho_{\widetilde{f}}(G_\Q)$
\end{itemize}

The hypothesis \eqref{item_irr} is vacuous when $N^-\neq1$. It also follows from \eqref{item_CR_3}(ii) if some prime divides $N^+$ exactly once, since the residual inertia image at such a prime contains a nontrivial unipotent of order $p$.

\begin{remark}\label{rk:EC_conditions}
For an elliptic curve $E/\QQ$ as in the introduction, we take $f$ to be the $p$-ordinary form specified in the introduction.
%The hypothesis \eqref{item_Heeg} requires the inert primes of bad reduction away from $p$ to be multiplicative, while all other bad primes away from $p$ split in $K$.
If $\epsilon^2\neq \mathbb 1$, this excludes the case where a quadratic twist gives good ordinary reduction. In particular, $\epsilon$ has order $3$, $4$ or $6$ by \cite[\S1.3]{delbourgo_Iwasawa_unstable}.
Writing $\epsilon=\omega^k$, condition \eqref{item_CR_1}(i) requires
\[
 2k<p-3,\qquad  \frac{p-1}{\gcd(p-1,2k+1)}>5.
\]
If $\epsilon$ has order $2$, \eqref{item_CR_1}(ii) requires that $p\ge7$ and $a_p(f)^2\not\equiv1\mod\varpi$.
Condition \eqref{item_CR_2} is equivalent to the absolute irreducibility of $E[p]$.
Let $\Delta_E$ denote the minimal discriminant of $E$.
At a multiplicative prime $\ell\neq p$, Tate uniformization identifies the ramification condition in \eqref{item_CR_3} with $p\nmid\ord_\ell(\Delta_E)$; at additive primes away from $p$, residual ramification is automatic since $p>3$.
Condition \eqref{item_CR_4} is the additional inertia-irreducibility condition at primes $\ell^2\mid N^+$ with $\ell\equiv1\pmod p$, and is vacuous when $N^+$ is square-free.
Finally, \eqref{item_irr} requires $p\mid\#\overline\rho_{E,p}(G_{\mathbf Q})$ when $N^-=1$; it is automatic if there is a multiplicative prime away from $p$ satisfying \eqref{item_CR_3}.
\end{remark}

\subsection{Selmer groups} \label{subsection_Selmer_groups}
\subsubsection{}
Let $V_{\widetilde{f}} := V_f \otimes \epsilon^{-1}$ be the self-dual twist of $V_f$. Fix a Galois-stable $\mathcal{O}$-lattice $T_f \subset V_{f}$, and set
$$ T_{\widetilde{f}}:= T_f \otimes \epsilon^{-1}, \qquad A_{\widetilde{f}} := V_{\widetilde{f}}/T_{\widetilde{f}}, \qquad T_n := T_{\widetilde{f}}/\varpi^n T_{\widetilde{f}}, \qquad A_n := A_{\widetilde{f}}[\varpi^n].$$
If we let $E = \Frac(\mathcal{O})$, then we have isomorphisms $T_n \cong A_n\simeq \Hom_\mathcal{O}(T_n, E/\mathcal{O}(1))$ as $G_{\mathbb{Q}}$-modules, where the first isomorphism is determined by our choice of uniformizer $\varpi$.

\subsubsection{\texorpdfstring{Local conditions at $\ell\nmid \Delta p$}{Local conditions at l}}
Let $\ell \neq p$ be a rational prime, $L/\mathbb{Q}$ finite, $w$ a prime of $L$ above $\ell$, and $M$ a discrete $G_L$-module. We set
$$H^1_{\rm{fin}}(L_w, M) := \ker \left(H^1(L_w, M) \lra H^1(I_w, M) \right), \qquad
H^1_{\rm{sing}}(L_w, M) := H^1(L_w, M)/H^1_{\rm{fin}}(L_w, M),$$ 
and write $\partial_w: H^1(L, M) \lra H^1_{\rm sing}(L_w, M)$ for the natural map induced by restriction composed with quotienting by $H^1_{\rm fin}$.

\subsubsection{\texorpdfstring{Local conditions at $p$}{Local conditions at p}}
By \eqref{item_ord}, we have an exact sequence of $\mathcal{O}[G_{\mathbb{Q}_p}]$-modules
\begin{equation}\label{eq:filtration}
    0 \lra F^+_p T_{\widetilde{f}} \lra T_{\widetilde{f}} \lra
F^-_p T_{\widetilde{f}} \lra 0,
\end{equation} 
where $F^\pm_p T_{\widetilde{f}}$ are free of rank one over
$\mathcal{O}$, and $G_{\mathbb{Q}_p}$ acts on $F^+_p T_{\widetilde{f}}$ by $\epsilon\,\chi_{\cyc}\,a_p^{-1}$ and on $F^-_p T_{\widetilde{f}}$ by $\epsilon^{-1}a_p$, where $a_p$ denotes the unramified character sending the arithmetic Frobenius to $a_p(f)$. We write $F^\pm_p M$ for the induced sub/quotients of $M \in \{T_n, A_n\}$, and define, for $L/K$ finite and $v \mid p$,
$$H^1_{\rm{ord}}(L_v, M) := \ker\left( H^1(L_v, M) \lra H^1(L_v, F^-_p M) \right).$$

\begin{lemma}\label{lemma_Howard_hyp_a}
Assume \eqref{item_Heeg}, \eqref{item_ord}, \eqref{item_CR_1} and \eqref{item_CR_2}. Then
$\overline\rho_{\widetilde f}|_{G_K}$ is absolutely
irreducible.
\end{lemma}

\begin{proof}
Suppose otherwise. Since twisting preserves absolute
irreducibility, Clifford theory gives
\[
 \overline\rho_{\widetilde f} \otimes\overline{\Fp}
 \cong \operatorname{Ind}_{G_K}^{G_{\QQ}}\psi
\]
for a character $\psi$.
Hence $\ker\overline\rho_{\widetilde f}\subseteq G_K$, so $K/\mathbf Q$ is unramified outside $Np$.
By \eqref{item_Heeg}, $(D_K,N)=1$, and therefore
\[
 K=\QQ\left(\sqrt{(-1)^{(p-1)/2}p}\right).
\]
Choose $\sigma\in I_p$ such that $u=\omega(\sigma)$ generates $\mathbf F_p^\times$.
Then $\sigma\notin G_K$, so the induced representation has trace zero at $\sigma$. On the other hand, \eqref{eq:filtration} gives
\[
 \operatorname{Tr}\overline\rho_{\widetilde f}(\sigma) =u^{k+1}+u^{-k} =u^{-k}(1+u^{2k+1})\neq0,
\]
since $u^{2k+1}$ has order greater than $5$ by \eqref{item_CR_1}. This is a contradiction.
\end{proof}

\begin{remark}\label{rk:cr2}
The proof of Lemma~\ref{lemma_Howard_hyp_a} shows that $\overline\rho_f$ remains absolutely irreducible over $\QQ\left(\sqrt{(-1)^{(p-1)/2}p}\right)$.
Thus, the irreducibility hypothesis in \cite{chidahsiehanticyclomainconjformodformscomposito} follows from \eqref{item_ord}, \eqref{item_CR_1} and \eqref{item_CR_2}.    
\end{remark}
\subsubsection{\texorpdfstring{Local conditions at $\ell \mid \Delta$}{Local conditions at l}}
Let $\ell$ be a prime such that $\ell \mid \Delta$, so either $\ell \mid N^-$ or $\ell$ is $n$-admissible. In either case (\textit{cf.} \cite{Carayol_representations_l_adique} when $\ell \mid \mid N^-$, Definition \ref{def_n-admissible_primes} when $\ell$ is $n$-admissible, and see also the discussion in \cite[\S1.2]{chidahsiehanticyclomainconjformodformscomposito}), we have an exact sequence of $\mathcal{O}[{G_{\Q_\ell}}]$-modules 
\begin{equation}
    0 \lra F_\ell^+ T_n \lra T_n \lra F_\ell^-T_n \lra 0, \label{eqn_Filtration_away_at_ell}    
\end{equation}
where $F^\pm_\ell T_n \cong \mathcal{O}_n$ as $\mathcal{O}$-modules and $G_{\Q_\ell}$ acts  on $F_\ell^+ T_n$ by $\tau_\ell \chi_\cyc$ and on $F_\ell^- T_n$ by $\tau_\ell$, where $\tau_\ell(\Frob_\ell) = \pm 1$ is an unramified character of $G_{\Q_\ell}$. As before, we write $F_\ell^\pm M$ for the induced sub/quotients of $M \in \{T_n, A_n \}$, and set 
$$H^1_{\rm{ord}}(L_w, M) := \ker\left( H^1(L_w, M) \lra H^1(L_w, F^-_\ell M) \right).$$
As explained in \cite[\S 2.2]{BertoliniDarmon2005}, when $\ell$ is $n$-admissible, the exact sequence \eqref{eqn_Filtration_away_at_ell} splits and we have a decomposition of $\mathcal{O}[G_{K_\ell}]$-modules
\begin{equation}
    T_n|_{G_{K_\ell}} = \mathcal{O}_n(1) \bigoplus \mathcal{O}_n. \label{eqn_decomp_of_T_at_ell}
\end{equation}

\subsubsection{Selmer groups}
For any prime $q$, we set 
\begin{align*}
    &H^1_\star(L_q, M):= \bigoplus_{w\mid q}H^1_\star(L_w,M), \quad \star \in \{\emptyset,  {\rm fin}, {\rm sing}, {\rm ord}\},
\end{align*}
where the sum is over primes $w$ of $L$ dividing $q$, and define $\partial_\ell$ appropriately.

\begin{defn}\label{defn_selmer}
    For $M \in \{T_n, A_n\}$ and $L/K$ finite, the Selmer group $\Sel_{\Delta}(L, M) \subset H^1(L,M)$ is defined by
    $$\Sel_{\Delta}(L, M):= \ker \left(H^1(L,M) \xrightarrow{(\oplus_{\ell \nmid p \Delta }\partial_\ell, \oplus_{\ell\mid p\Delta} \res_\ell)} \bigoplus_{\ell\nmid p\Delta }H^1_{\rm sing}(L_\ell, M) \oplus \bigoplus_{\ell \mid p\Delta}\frac{H^1(L_\ell, M)}{H^1_{\rm ord}(L_\ell, M)} \right).$$
    In particular, the Selmer group consists of the classes $s$ with $\partial_\ell(s) = 0$ for all $\ell \nmid p \Delta $ and $\res_\ell(s) \in H^1_{\rm ord}(L_\ell, M)$  for all $\ell \mid p\Delta$. We also set 
    \begin{align*}
        &\Sel_\Delta(K_{\rm{ac}}, A_n):= \varinjlim_m \Sel_\Delta(K_m, A_n), &&\Sel_\Delta(K_\ac, T_n):= \varprojlim_m \Sel_\Delta(K_m, T_n), \\
        &\Sel_{N^-}(K_{\rm{ac}}, A_{\widetilde{f}}):= \varinjlim_n \Sel_{N^-}(K_\ac, A_n), &&\Sel_{N^-}(K_\ac, T_{\widetilde f}):= \varprojlim_n \Sel_{N^-}(K_\ac, T_n),\\
                &\Sel_{N^-}(K_m, A_{\widetilde{f}}):= \varinjlim_n \Sel_{N^-}(K_m, A_n), &&\Sel_{N^-}(K_m, T_{\widetilde{f}}):= \varprojlim_n \Sel_{N^-}(K_m, T_n).
    \end{align*}
    %    The groups $\Sel_\Delta(K_\ac, A_n)$, $\Sel_{N^-}(K_\ac $, $ A_{\widetilde{f}})$ (resp. $\Sel_\Delta(K_\ac, T_n)$, $\Sel_{N^-}(K_\ac, T_{\widetilde{f}})$) are discrete (resp. compact) $\Lambda(\Gamma_\ac)$-modules.
\end{defn}

\subsection{Local cohomology}
Here we gather some useful local cohomology results. 
\subsubsection{\texorpdfstring{Results at $p$}{Results at p}}

\begin{lemma}\label{lemma_vanishing_cohomology_at_p}
    Let $v$ be a prime of $K_m$ above $p$ and $K_{m,v}$ be the completion at $v$. Then
        \begin{itemize}
            \item[i)] $H^0(K_{m,v}, F_p^-A_{\widetilde{f}})=H^2(K_{m,v}, F_p^+T_{\widetilde{f}})=0$.
            \item[ii)] $H^0(K_{m,v}, F_p^+A_{\widetilde{f}}) = H^2(K_{m,v}, F_p^-T_{\widetilde{f}})=0 $.
        \end{itemize}
\end{lemma}

\begin{proof}
    (i) Suppose first that $\epsilon^2 \neq \mathbb 1$ or $p\nmid D_K$. Let $I_p$ denote the inertia subgroup of $G_{\Q_p}$, and $I_v \subset I_p$ denote the inertia subgroup of $G_{K_{m,v}}$. It suffices to show that $F_p^-T_n$ has no $I_v$-invariants. Since $I_v$ acts on $F_p^-T_n$ through $\epsilon^{-1}$, and $\epsilon$ has order prime to $p$, we are reduced to showing that $\epsilon|_{I_v}$ is non-trivial. Let $d$ denote the order of $\epsilon(I_p)$. By the isomorphism 
    $$I_p/ I_v\cdot \ker(\epsilon) \cong \mu_d/ \epsilon(I_v),$$
    we are reduced to showing $[I_p: I_v \cdot \ker(\epsilon)]< d$. When $\epsilon^2 \neq \mathbb 1$, this follows since $d>2$ and the prime-to-$p$ part of $[I_p:I_v]$ divides $2$. If $p\nmid D_K$, then $[I_p:I_v] $ is a power of $p$, so $\epsilon|_{I_v}$ remains non-trivial.
    
    Suppose $p\mid D_K$ and \eqref{item_CR_1}(ii) holds. Let $w$ be the place of $K$ lying below $v$. Then $F_p^-T_1$ is unramified over $K_w$, with Frobenius eigenvalue $\pm\overline{a_p(f)}\ne 1$. Since $K_{m,v}/K_w$ has $p$-power degree, the argument of \cite[Lemma 1.8]{chidahsiehanticyclomainconjformodformscomposito} gives $H^0(K_{m,v}, F_p^-A_{\widetilde{f}})=0$ and the $H^2=0$ statement follows by local duality.

    For (ii), the action of $I_p$ on $F_p^+T_1$ is given by
    $\epsilon\chi_{\mathrm{cyc}}
    =\omega^{k+1}$. {By \eqref{item_CR_1}, $\omega^{k+1}$ has order $e>2$ and $p\nmid e$}. If $\omega^{k+1}|_{I_v}$ is trivial, then $e\, |\, [I_p:I_v]\, |\, 2p^m$, which is impossible. Hence, $H^0(K_{m,v}, F_p^+T_1)=0$. It then follows that $H^0(K_{m,v}, F_p^+A_{\widetilde{f}})=0$, and the result for $H^2$  follows by local duality.
\end{proof}

\begin{proposition}\label{prop_orthogonal_at_p}
    $H^1_\ord (K_{\ac,p}, A_n )$ and $H^1_\ord(K_{\ac, p}, T_n)$ are orthogonal complements under the local Tate pairing. 
\end{proposition}
\begin{proof}
    This follows as in the proof of \cite[Proposition 1.7 (1)]{chidahsiehanticyclomainconjformodformscomposito}, using the previous lemma.
\end{proof}

\subsubsection{\texorpdfstring{Results at $\ell \neq p$}{Results at }}

\begin{theorem} \label{theorem_Howard_hyp_e}
    Let $s\leq n $ be positive integers and let $\kappa \in H^1(K, A_s)$ be a non-zero element. Then there exist infinitely many $n$-admissible primes $\ell$ satisfying $\partial_\ell(\kappa) = 0$, and the map 
    $$\res_\ell: \langle \kappa\rangle \lra H^1_{\rm fin}(K_\ell, A_s)$$
    is injective, where $\langle \kappa \rangle \subset H^1(K, A_s)$ is the $\mathcal{O}$-submodule generated by $\kappa$.
\end{theorem}

\begin{proof}
    Let us set 
    $$\rho_n := \rho_{\widetilde{f}} \mod \varpi^n, \qquad F_n := \overline \Q^{\ker \rho_n}.$$
    If $K$ and $F_n$ are linearly disjoint for all $n$, then the theorem follows exactly as in the proof of \cite[Theorem 6.3]{chidahsiehanticyclomainconjformodformscomposito}. Now suppose that $K \subset F_n$. Since $(D_K,N)=1$ and $\rho_n$ is unramified outside $pN$, we must have that $D_K=p$, and hence we have the equality
    $$ \chi_K(\sigma) = \left(\frac{\omega(\sigma)}{p}\right) = \left(\frac{\det \rho_n(\sigma)}{p}\right),$$
    where $\chi_K$ denotes the quadratic Dirichlet character of conductor $p$ attached to $K$.
    We indicate the necessary changes to the argument in \textit{loc. cit.} so that the theorem also holds in this case, using the notation of \textit{loc. cit.} freely. Following \cite[proof of Theorem~6.3]{chidahsiehanticyclomainconjformodformscomposito}, we may assume that $s=1$ and $\tau\kappa=\delta\kappa$, where $\delta\in\{\pm1\}$ and $\tau$ is the complex conjugation. We have a decomposition $\Gal(M_\kappa/ \Q) \cong A_\kappa \rtimes \Gal(F_n/\Q)$, where in our case, the action of $\Gal(F_n/\Q)$ on $A_\kappa$ is given by $\sigma \cdot v = \delta^{j(\sigma)}\rho_1(\sigma) \cdot v$, where $\delta =\pm 1$, $v \in A_\kappa$, $(-1)^{j(\sigma)} = \chi_K(\sigma)$. Note that in the situation considered in \textit{loc. cit}, $j(\sigma)=j$ is independent of $\sigma$. We may choose $\sigma \in \Gal(F_n/\Q)$ such that $\rho_n(\sigma)$ has eigenvalues $\delta$, $\lambda$, with $\pm 1 \neq \lambda \in (\Z/p^n\Z)^\times$, and $\lambda$ has order prime to $p$. Moreover, we can also ensure that 
    \begin{equation}
        \chi_K(\sigma) = -1. \label{eqn_legendre=-1} \notag
    \end{equation} Indeed, $\rho_n(\sigma)$ is determined by an element $h \in \overline \rho_{\widetilde f}(\Gal(F_n/\Q))$, which, as is explained in the proof of \cite[Lemma 6.2]{chidahsiehanticyclomainconjformodformscomposito}, can be taken to be either
    $$h = \begin{cases}
        \pm \overline \rho_{\widetilde f}(\Frob_\ell), &\text{if there exists } \ell\mid N^- : p \nmid (\ell^2 -1),\\
        \text{conjugate to some } \gamma \in \GL_2(\Fp): \det(\gamma) \in \Fp^\times \setminus \left((\Fp^\times)^2\cup \{-1\}\right) , &\text{if }p\mid \#\overline \rho_{\widetilde{f}}(G_\Q).
    \end{cases}$$
    In either of these cases, we have that $\chi_K(\sigma) = \left(\frac{\det(h)}{p}\right)=-1$ and the result follows by the same Chebotarev density argument in \textit{loc. cit.}
\end{proof}

We set
\begin{align*}
    &H^1(K_{\ac}, T_n) := \varprojlim_m H^1(K_m, T_n),
    &H^1(K_{\ac}, A_n) := \varinjlim_m H^1(K_m, A_n).
\end{align*}
If $M \in \{T_n, A_n\}$, $\ell$ is a prime of $K$, and $L$ is a finite extension of $K$, then we write 
$$H^1(L_{\ell}, M) := \bigoplus_{w \mid \ell}H^1(L_w, M), $$
where the sum is over primes of $L$ dividing $\ell$.

\begin{lemma}\label{lemma_cohomology_decomp_at_admissible_ell}
    Let $\ell$ be an $n$-admissible prime. Then 
    we have decompositions
    \begin{align*}
        &H^1(K_{\ac, \ell}, T_n) = H^1_{\rm fin}(K_{\ac, \ell}, T_n) \bigoplus H^1_\ord(K_{\ac, \ell}, T_n), \\
        &H^1(K_{\ac, \ell}, A_n) = H^1_{\rm fin}(K_{\ac, \ell}, A_n) \bigoplus H^1_\ord(K_{\ac, \ell}, A_n).
    \end{align*}
    Moreover, for $\star \in \{{\rm ord}, {\rm fin}\}$, we have that $H^1_\star(K_{\ac, \ell}, T_n)$ is a maximal isotropic submodule under the Tate pairing, and an isomorphism of $\mathcal{O}_n[[\Gamma_\ac]]$-modules 
    \begin{equation}
        H^1_\star(K_{\ac, \ell}, T_n) \cong \mathcal{O}_n[[\Gamma_\ac]]. \label{eqn_isomorphism_H^1=Lambda_n}   
    \end{equation}
    
\end{lemma}

\begin{proof}
    This follows from combining \eqref{eqn_decomp_of_T_at_ell} with Lemmas 2.2.1 and 3.1.2 of \cite{howard06} which are stated in the context of elliptic curves, but whose proofs clearly generalize to our case.
\end{proof}

\begin{lemma} \label{lemma_orthogonal_complements_at_ell_dividing_N}
    Let $\ell \mid N^-$. Then $H^1_\ord(K_{\ac,\ell}, T_n)$ and $H^1_\ord(K_{\ac,\ell}, A_n)$ are orthogonal complements under the local Tate pairing. 
\end{lemma}

\begin{proof}
    This follows exactly as in the case $\ell \mid N^-$ in the proof of \cite[Proposition 1.7 (1)]{chidahsiehanticyclomainconjformodformscomposito}.
\end{proof}

\section{\texorpdfstring{CM points, $p$-adic $L$-functions, and explicit reciprocity laws}{CM points, p-adic L-functions and explicit reciprocity laws}}

\subsection{CM points}
Let $\Delta \in \Z$ be as in \S \ref{section_Local_cohomology_and_Selmer_groups} and fix $B_\Delta$ to be the quaternion algebra over $\Q$ of discriminant $\Delta$.

\subsubsection{CM points on Quaternionic Shimura varieties}

\begin{lemma}\label{lemma_optimal_embeddings}
    Let $R^0 \subset R$ be the Eichler order of level $pN^+$ whose completion at $p$ is the standard Eichler order of level $p\Z_p$. There exists an embedding $\jmath : K \hookrightarrow B_\Delta$ and, for each $m \geq 0$, an element $\varsigma_m \in \widehat{B}_\Delta^\times$ with
    $\varsigma_{m,v} = 1$ for $v \neq p$ and $\varsigma_{m,p} =
    \iota^{-1}_p(\gamma\cdot\left(\begin{smallmatrix} p^m & 0 \\ 0 & 1
    \end{smallmatrix}\right))$ for a fixed $\gamma \in \GL_2(\mathbb{Z}_p)$, such that:
    \item[i)] $\jmath(K) \cap R = \jmath(\mathcal{O}_K)$,
    \item[ii)] For $m\geq 1$, we have the equality of sets
    $$\jmath(K) \cap \varsigma_m \widehat{R}^0\, \varsigma_m^{-1}
     =\jmath(K) \cap \varsigma_m \widehat{R}\, \varsigma_m^{-1} = \jmath(\mathcal O_K[p^m]).$$
    In particular, the order $B_\Delta \cap \varsigma_m \widehat R^0 \varsigma_m^{-1}$ is `admissible for $(\pi_f, p^m, \jmath)$', in the sense of \cite[Definition 1.4]{kohen_pacetti_anticyc_l-function}
    \item[iii)] For all $m\geq n\geq 1$, we have the containment 
    $$\varsigma_m^{-1}\,
\widehat{\jmath}(\widehat{\mathcal{O}_K[p^m]}^\times)\, \varsigma_m \subset
\su_n.$$ 
\end{lemma}
\begin{proof}
    Part (i) is standard under our hypotheses (see \cite{Gross_88, Vigneras_Quaternions}), and (ii) is a slight restatement of \cite[Lemmas 1.6, 1.7]{kohen_pacetti_anticyc_l-function}, however, we give a proof of (ii) as it will be helpful in the proof of (iii).
    
    For (ii), it suffices to work only at $p$, since $\mathcal{O}_K[p^m] \otimes \Z_q \simeq \mathcal{O}_K \otimes \Z_q$ for $q \nmid p$. Let $\omega_K \in K$ be such that $\mathcal{O}_K = \Z[\omega_K]$. We fix a $\Z$-basis $\{\omega_K, 1\}$ of $\mathcal{O}_K$, and to ease notation, we set $B=B_\Delta$, $\mathcal{O}_p[p^m]:= \mathcal{O}_K[p^m] \otimes \Z_p$ and $\mathcal{O}_p := \mathcal{O}_K \otimes \Z_p$. Then the isomorphism $K_p:= K \otimes \Q_p \cong \Q_p \cdot \omega_K \bigoplus \Q_p  $ determines an embedding $$\rho_p: K_p \hookrightarrow \End_{\Q_p}(K_p) \cong \iota^{-1}_p(M_2(\Q_p)) = B_p,$$
    where the first arrow maps $k$ to `multiplication by $k$', and we have
    \begin{equation}
        \rho_p(K_p) \cap R_p = \rho_p(\mathcal{O}_p) \label{eqn_optimal_embedding_p}.
    \end{equation}
    Writing $\delta_m = \begin{psmallmatrix}
        p^m & 0\\
        0 & 1
    \end{psmallmatrix}$, a quick computation shows that for $m\geq 0$, our choice of basis gives an isomorphism
    \begin{equation}
        \delta_m M_2(\Z_p) \delta_m^{-1} \cong \End_{\Z_p}(\mathcal{O}_p[p^m]), \label{eqn_Endomorphism_ring_identity}
    \end{equation}
    and hence $$\rho_p(K_p) \cap \iota^{-1}_p(\delta_m)R_p\iota^{-1}_p(\delta_m)^{-1} = \rho_p(\mathcal{O}_p[p^m]).$$
    For $m\geq 1$ we want to show the same equality, but with $R^0_p$ instead of $R_p$. Since $R^0_p \subset R_p$, it suffices to show the inclusion 
    $$\rho_p(\mathcal{O}_p[p^m]) \stackrel{?}\subset \rho_p(K_p) \cap \iota^{-1}_p(\delta_m)R^0_p\iota^{-1}_p(\delta_m)^{-1}.$$Recall the identification \eqref{eqn_Endomorphism_ring_identity} given by our fixed basis $\{\omega_K, 1\}$ and note that the change of basis matrix $\{\omega_K, 1\} \mapsto \{p^m \omega_K, 1\}$ is $\delta_m$. Hence, this new choice of basis corresponds to the embedding $\rho_m := \iota^{-1}_p(\delta_m)^{-1} \rho_p \iota^{-1}_p(\delta_m) $, and we have 
    $$\rho_p(x)\in \iota_p^{-1}(\delta_m) R^0_p\iota_p^{-1}(\delta_m)^{-1} \Longleftrightarrow \rho_m(x) \in R^0_p.$$
    For $x\in \mathcal{O}_p[p^m]$, there exist unique $r,s \in \Z_p$ such that the matrix associated to $x = r\cdot p^m\omega_K + s$ (viewed as an element of $ \End_{\Z_p}(\mathcal{O}_p[p^m]) $ using $\{p^m\omega_K,1\}$-coordinates) is given by 
    \begin{equation}
    \begin{pmatrix}
        s+ rp^m \cdot \Tr_{K/\Q}(\omega_K) & r\\
        -rp^{2m}\cdot \rm{N}_{K/\mathbb{Q}}(\omega_K) & s
    \end{pmatrix} \equiv
    \begin{pmatrix}
        s & r\\
        0  & s
    \end{pmatrix} \mod p^m, \label{eqn_matrix_for_x_in_new_coords}
    \end{equation}
    which belongs to $\iota_p(R^0_p)$, as required.
    The result now follows on noting that \eqref{eqn_optimal_embedding_p} implies the existence of an element $\gamma \in \GL_2(\Z_p)$ such that $\gamma^{-1} \jmath_p \gamma = \rho_p$. 

    (iii) Again, it suffices to show the containment at $p$, and the result is now immediate by the previous discussion, since every element in $\varsigma_{m,p}^{-1} \,  \jmath_p(\mathcal{O}_p[p^m]) \, \varsigma_{m,p}$ corresponds to a matrix of the form given in \eqref{eqn_matrix_for_x_in_new_coords}.
\end{proof}

Henceforth, for a given $\Delta$, we fix a tuple $(\jmath, (\varsigma_m)_m)$ satisfying the conditions of the previous lemma in a compatible way for varying $\Delta$. If $B_\Delta$ is indefinite, then we set $z_0 \in \mathcal{H}^+$ to be the unique element in the upper half-plane fixed by the action of $\jmath(K^\times)$ through $\iota_\infty$.

\begin{defn}(CM points)\label{def_cm-points}
Let $U \subset \widehat{R}^\times$ be an open compact subgroup and
$m \geq 0$. For $a \in \widehat{K}^\times$, the \emph{CM point $x_m(a)$ of
conductor $p^m$} is defined by
$$ x_m(a) := \begin{cases}
[\widehat{\jmath}(a)\,\varsigma_m]_{U} \in X(U)(\C), &
\text{$B_\Delta$ definite},\\
[(z_0, {\widehat{\jmath}}(a)\,\varsigma_m)]_{U} \in X(U)(\C), & \text{$B_\Delta$ indefinite}.
\end{cases} $$
In the definite case, we set $\mathrm{CM}(U)=\widehat{\jmath}(K^\times)\backslash\widehat B^\times_\Delta/U$ and $\pr_U:\mathrm{CM}(U)\to B_\Delta^\times\backslash\widehat B^\times_\Delta/U$.
\end{defn}

Suppose $B_\Delta$ is indefinite, and let $\rec_K: K^\times \backslash \widehat K^\times \xrightarrow{\sim} \Gal(K^{\rm ab} /K )$ denote the geometrically normalized Artin map. Then Shimura's reciprocity law states
\begin{equation}
    \rec_K(a)[z,b]_U = [z, \widehat \jmath(a)b ]_U \label{eqn_Shimura_reciprocity}.
\end{equation}
If $B_\Delta$ is definite, $G_K$ acts on $\mathrm{CM}(U)$ via  
\begin{equation}
    \rec_K(a)[b]_U = [ \widehat \jmath(a)b ]_U \label{eqn_definite_Galois_action}.
\end{equation}
This induces an action on $x_m(a)$ via $\pr_U$.
\begin{remark}\label{remark_fields_of_definition_of_CM_points}
    Let $U = \su_n$ in Definition \ref{def_cm-points}. From Lemma \ref{lemma_optimal_embeddings} (iii), we see that the CM point attached to $a \in \widehat K^\times$ depends only on the class $[a]\in K^\times \backslash \widehat K^\times / \widehat{\mathcal{O}}_K[p^m]^\times \simeq \Gal(H_m/K)$ for $m\geq n\ge 1$. Combining this with \eqref{eqn_Shimura_reciprocity} and \eqref{eqn_definite_Galois_action} yields $x_m(a)\in X(U)(H_m)$ in the indefinite case. In the definite case, the CM parametrization $a\mapsto x_m(a)$ factors through $\Gal(H_m/K)$.
\end{remark}

Recall the double coset decomposition \eqref{eqn_double_coset_decomp}. This defines a Hecke action on ${\rm Div}(X(\su_n))$, given by 
\begin{align*}
    &[\su_n x \su_n]\cdot [b]_{\su_n} = \sum_{i\in I}[bx_i]_{\su_n},
    && B_\Delta \text{ definite,} \\
    &[\su_n x \su_n]\cdot [z,b]_{\su_n} = \sum_{i\in I}[z,bx_i]_{\su_n},
    && B_\Delta \text{ indefinite.}
\end{align*}

\begin{proposition}[Distribution relation]\label{prop_distribution_relation}
For $m\ge n\ge1$, we have the following equality in $\mathrm{Div}(X(\su_n))$,
    $$\sum_{\sigma \in \mathrm{Gal}(H_{m+1}/H_m)} x_{m+1}(a)^\sigma
    = U_p \cdot x_m(a).$$ 
\end{proposition}

\begin{proof}
    This is essentially \cite[Proposition 3.4]{longovigni_11}; although their setting is slightly different as they use $\Gamma_1(p)$-type level subgroups, the computation is the same.
\end{proof}

\subsubsection{Heegner classes}\label{s:heegnerclasses}
Suppose that $\omega(\Delta)$ is even, so that $B_\Delta$ is indefinite. Recall that $\mathcal{O}$ denotes the ring of integers of $\Q_p((a_n(f))_n)$.
Fixing an $n$-admissible form $\mathcal{D}= (\Delta,g)$ (\textit{cf.} Definition \ref{def_n-admissible_forms}), we write $I_g = \ker(g: \T_{\Delta} \twoheadrightarrow \mathcal{O}_n)$.  

Suppose $m\geq n\geq 1$ and choose a prime $q_0 \nmid p \Delta N^+$ such that $1+ q_0 - \epsilon^{-1}(q_0)\lambda_{q_0}(f) \in \mathcal{O}^\times$ and $\left(\frac{q_0}{p}\right)=1$.
Such $q_0$ exists by Remark~\ref{rk:cr2} and Chebotarev.
Setting $P_m = x_m(1)$, we define a class 
\begin{equation}
    D_m := \sum_{\sigma \in \Gal(H_m/ K_m)} (1+q_0 -T_{q_0})\cdot P_m^\sigma \otimes  (1+q_0 -\epsilon^{-1}(q_0)\lambda_{q_0}(f))^{-1} \in J(\su_n)(K_m) \otimes \mathcal{O},\label{eq:CM}
\end{equation}
where $J(\su_n)$ denotes the Jacobian variety of $X(\su_n)$.
By the description of the connected components given in \S\ref{s:components}, the Hecke operator $T_{q_0}$ (defined as in \S\ref{S:hecke}) preserves each component. Furthermore, the degree of the divisor $T_{q_0}D$ is equal to $(q_0+1)\deg D$ for any divisor $D$. Thus, $D_m$ has degree zero on each component. 

If $T_p(J(\su_n))$ denotes the $p$-adic Tate module of $J(\su_n)$, then we define the cohomology classes 
$$\kappa_{\mathcal{D},m} = a_p(f)^{-m}\operatorname{Kum}(D_m) \mod I_g \in H^1(K_m, (T_p(J(\su_n))\otimes \mathcal{O})/I_g) \simeq H^1(K_m, T_n),$$
where $\operatorname{Kum}:J(\su_n)(K_m) \otimes\mathcal{O} \lra H^1(K_m, T_p(J(\su_n))\otimes\mathcal{O})$ is the Kummer map and the final isomorphism follows from \cite[Corollary 4.4]{chidahsiehanticyclomainconjformodformscomposito}, \cite[Theorem 5.17]{BertoliniDarmon2005}, \cite[Proposition 4.4]{pollack-weston11}.
The classes $\kappa_{\mathcal{D},m}$ are compatible under corestriction maps as $m$ varies by Proposition \ref{prop_distribution_relation}, and are independent of the choice of $q_0$.

\begin{defn}[Anticyclotomic Heegner class]
    Let $\mathcal{D} = (\Delta,g)$ be an $n$-admissible form. The Heegner class $\kappa_\mathcal{D}$ attached to $\mathcal{D}$ is defined as 
    $$\kappa_\mathcal{D} := \varprojlim_m \kappa_{\mathcal{D},m} \in H^1(K_\ac, T_n).$$
\end{defn}

Suppose now that there exists an $n$-admissible prime $\ell \mid \Delta$, and set $\Delta^{(\ell)} := \frac{\Delta}{\ell}$, so that $\omega(\Delta^{(\ell)})$ is odd. Suppose $ \mathcal{D}^{(\ell)} = (\Delta^{(\ell)}, g^{(\ell)})$ is an $n$-admissible form. In \cite[\S 4.3.2]{chidahsiehanticyclomainconjformodformscomposito}, the authors attach a Heegner class $\kappa_{\mathcal{D^{(\ell)}}}(\ell)$ to the data $(\mathcal{D}^{(\ell)}, \ell)$. The following lemma states that our Heegner classes coincide with theirs.

\begin{lemma} \label{lemma_comparing_Heegner_classes}
    With notation as above, for all $n$-admissible primes $\ell \mid \Delta$, we have the equality
    $$\kappa_\mathcal{D} = \kappa_{\mathcal{D}^{(\ell)}}(\ell).$$
\end{lemma}

\begin{proof}
    We freely use the notation in \cite{chidahsiehanticyclomainconjformodformscomposito}. Let $\varphi_{B_{\Delta^{(\ell)}}, B_{\Delta}}$ be the map defined in \S 3.1 of \textit{op. cit.} Let $\mathscr{U}_n^{(\ell)}\subset \widehat B^\times_{\Delta^{(\ell)}}$ be the open compact subgroup defined as in \S \ref{subsub_QuaternionPrelims}, but without the component at $\ell$. Then it is clear that the open compact subgroup $\varphi_{B_{\Delta^{(\ell)}}, B_{\Delta}}(\mathscr{U}_n^{(\ell)}) \cdot \mathcal{O}^\times_{B_{\Delta, \ell}}$ coincides with $\mathscr{U}_n \subset \widehat B^\times_{\Delta}$, so the corresponding Shimura curves and CM points coincide. Moreover, the ideals $I_{g^{(\ell)}}^{[\ell]}$ and $I_g$ must be equal since the corresponding homomorphisms have the same eigensystem (where $I_g^{[\ell]}$ denotes the ideal defined in \S 4.2 of \textit{op. cit.}), and the result now follows from the definitions of $\kappa_\mathcal{D}$ and $\kappa_{\mathcal{D}^{(\ell)}}(\ell)$ (after possibly rescaling the isomorphism $T_p(J_n^{[\ell]})/I_{g^{[\ell]}} \simeq T_n$ by an element of $\mathcal{O}_n^\times$).
\end{proof}

\begin{proposition}
     Let $\mathcal{D} = (\Delta,g)$ be an $n$-admissible form. Then 
     $$\kappa_\mathcal{D} \in \Sel_{\Delta}(K_\ac, T_n).$$
\end{proposition}
\begin{proof}
    The proof of \cite[Proposition 4.7]{chidahsiehanticyclomainconjformodformscomposito} can be generalized to our current setting. The local conditions away from $p$ have been verified in \textit{loc. cit.} At $p$, the argument in the proof of Lemma~\ref{lemma_vanishing_cohomology_at_p}(i) gives us the corresponding vanishing condition in the proof therein. 
\end{proof}

\subsection{\texorpdfstring{The $p$-adic $L$-function}{The p-adic L-function}}\label{sec:KP}
With the same conventions as in \S \ref{section_Local_cohomology_and_Selmer_groups}, in this section we recall the anticyclotomic $p$-adic $L$-function constructed in \cite{kohen_pacetti_anticyc_l-function}, which is attached to a $p$-ordinary weight-two eigenform of level $\Gamma_0(N)\cap \Gamma_1(p)$. Although the authors only consider newforms attached to potentially ordinary elliptic curves, the construction generalizes to our setting, as we explain below.

\subsubsection{}\label{S:padicL}
Assume that $\omega(\Delta)$ is odd. Let $g \in S_2(\Gamma_0(\Delta N^+) \cap \Gamma_1(p), \psi^{-2})$ be a $p$-ordinary eigenform that is new away from $p$, with $U_p$-eigenvalue $\alpha \in \mathcal{O}^\times$ and with $\psi \neq \mathbb 1$. When $\psi^2=\mathbb{1}$,  let $g^\circ$ be the newform of level $\Delta N^+$ whose ordinary $p$-stabilization is $g$. For an anticyclotomic Hecke character $\chi$ of $K$, with $p$-power conductor, we write
$$\widetilde{\chi} := \chi \cdot (\psi \circ \mathrm{N}_{K/\mathbb{Q}}),$$
and, following \cite[Definition 1.2]{kohen_pacetti_anticyc_l-function}, we say $\chi$ is \emph{good} if $p$ divides the conductor of $\widetilde{\chi}$. As in \textit{op. cit.}, one has the following identity of $L$-functions
\begin{equation}\label{eqn_L-twist}
L(\pi_{g \otimes \psi}, \chi, s)
 = L(\pi_g, \widetilde{\chi}, s).
\end{equation}
By Lemma \ref{lemma_optimal_embeddings} (ii), there exists an embedding $\jmath:K \hookrightarrow B_\Delta$ such that for all $m\geq 1$, the order $B_\Delta \cap \varsigma_m \, \widehat R^0 \, \varsigma_m^{-1}$ is admissible for $(\pi_g, p^m , \jmath)$, hence we use the construction of \S1--3 (respectively \S4) of \textit{op. cit.} when $\psi^2\ne\mathbb{1}$ (respectively $\psi^2=\mathbb{1}$), which we spell out in the next subsubsection, after the following preliminaries. 
In \cite{kohen_pacetti_anticyc_l-function}, the authors consider the open compact subgroup of $\widehat B_\Delta^\times$ defined by 
$$\su^{\rm{KP}} := \{ x \in \widehat{R}^\times :
\iota_p(x_p) \equiv \begin{psmallmatrix} 
* & * \\
0 & a
\end{psmallmatrix} \mod p,\ \psi^{-2}(a) = 1 \}.$$

When $\psi^2\ne\mathbb{1}$,  let $g_{B_\Delta} : B_\Delta^\times \backslash \widehat{B}_\Delta^\times \lra \mathcal{O}$ be the Jacquet--Langlands transfer of $g$, as in \cite[Theorem 1.9]{kohen_pacetti_anticyc_l-function}, so that $g_{B_\Delta}$ is $\su^{\rm{KP}}$-invariant. When $\psi^2=\mathbb{1}$, we use the ordinary stabilization of the Jacquet--Langlands transfer of $g^\circ$ instead.
In both cases, let $F_\psi \subset \mathbb{Q}(\mu_p)$ be the field cut out by $\ker\psi$. We set 
$$
\widetilde{H}_m := H_m F_\psi, \qquad \widetilde{G}_m := \Gal(\widetilde{H}_m/K), \qquad \widetilde{\Gamma}_\ac := \varprojlim_m \widetilde{G}_m. $$
%\red{Yep. I have relabeled the earlier one as $\overline \Gamma$. In fact I think we can get rid of the twisting discussion in \S2 since it isn't very relevant to the rest of the paper. I wrote it in the very beginning when I thought that twisting the BDP $L$-function might be an approach.} \green{I see, that makes sense. Yes, perhaps getting rid of that would be a good idea if it is not used in the rest of the paper.}
\subsubsection{}\label{sec:p-adic-L}
We continue to assume that $\omega(\Delta)$ is odd. In \cite[\S 3 and \S4]{kohen_pacetti_anticyc_l-function}, the authors define theta elements by first setting $\zeta_m^\sigma := \alpha^{-m}\, g_{B_\Delta}(x_m(a_\sigma))$ for $\sigma = {\rm rec}_K(a_\sigma)|_{\widetilde{H}_m} \in \widetilde{G}_m$, and then defining
\begin{equation}\label{eqn_KP-theta}
    \theta_{g, m}^{\rm KP} := \sum_{\sigma \in \widetilde{G}_m} \psi(\sigma)\, \zeta_m^\sigma\cdot \sigma \in \mathcal{O}[\widetilde{G}_m], \qquad \theta_g^{\rm KP} := \varprojlim_m \theta_{g, m}^{\rm KP} \in \Lambda( \widetilde{\Gamma}_\ac),
\end{equation}
where $\psi$ is regarded as a character of
$\rm{Gal}(\overline{\mathbb{Q}}/\mathbb{Q})$ restricted to
$\widetilde{G}_m$. The $p$-adic $L$-function $\mathscr{L}_g $ is then given by the image of
\[ \theta_g^{\rm KP} \cdot
(\theta_g^{\rm KP})^{*} \in \Lambda( \widetilde{\Gamma}_\ac )\to\Lambda(\Gamma_\ac),\]
where $*$ denotes the involution $\sigma \mapsto \sigma^{-1}$ for $\sigma \in \widetilde{\Gamma}_\ac$. We then have the following theorem.

\begin{theorem}\label{theorem_p-adic_l-function_interpolation}
    The $p$-adic $L$-function $\mathscr{L}_g$ satisfies the following interpolation property: there is a complex period $\Omega_g \in \mathbb{C}^\times$ such that for
    every good anticyclotomic character $\chi$ with
    ${\rm cond}(\widetilde{\chi}) = p^m$, $m \geq 1$,
    $$\chi(\mathscr{L}_g) \doteq \frac{p^{m}}{\alpha^{2m}} \cdot
    \frac{L(g \otimes \psi/K, \chi, 1)}{\Omega_g} \cdot \sqrt{D_K},$$
    where $\doteq$ denotes equality up to an explicit $p$-adic unit and $L(g \otimes \psi/K, \chi, s) = L(\pi_{g \otimes \psi}, \chi, s)$ (\textit{cf.} \eqref{eqn_L-twist})

    If $\psi^2=\mathbb{1}$ and $\chi$ is a finite-order character of $\Gamma_\ac$ such that $\widetilde\chi$ is unramified at primes above $p$, then $p\cO_K=\fp^2$ and
    \[
    \chi(\mathscr{L}_g)\doteq\frac{(1-\alpha^{-1}\widetilde\chi(\fp))^2}{L_p(\pi_{g^\circ},\widetilde\chi,1)}\cdot\frac
    {L(g^\circ\otimes\psi/K,\chi,1)}{\Omega_g}\cdot\sqrt{D_K},
    \]
    where $L_p(\pi_{g^\circ},\widetilde\chi,1)$ denotes the product of the local Euler factors above $p$.
\end{theorem}
\begin{proof}
    This is \cite[Theorems 3.6 and 4.1]{kohen_pacetti_anticyc_l-function}. The proof of the interpolation formula rests on a combination of the explicit Waldspurger formula of Cai--Shu--Tian \cite[Theorem 1.8]{caishutian} for the pair $(\pi_g, \widetilde{\chi})$, which holds for a form $g$ of weight $2$ satisfying $\omega_{\pi_g}\cdot\widetilde{\chi}|_{\mathbb{A}^\times} = 1$, where $\omega_{\pi_g}$ denotes the central character, together with \cite[Proposition 3.5]{kohen_pacetti_anticyc_l-function}, which holds formally in our case as well. Under \eqref{item_Heeg} the set $\Sigma_D$ of \cite[Theorem 3.6]{kohen_pacetti_anticyc_l-function} is empty, and for $m \geq 1$, the unit factors $u_m = 1$, which yields the displayed formula.

    When $\psi^2=\mathbb{1}$, the construction of the $p$-adic $L$-function corresponds to the one for $g^\circ$ given in \cite[\S4]{kohen_pacetti_anticyc_l-function}. The argument of \cite[Theorem~4.1]{kohen_pacetti_anticyc_l-function} gives the second formula. The factor arising from $\cO_K^\times$ therein is a $p$-adic unit in our current setting.
\end{proof}

\begin{defn}\label{defn_KP_theta_classes}
    We set $\theta_{g, m} \in \mathcal{O}[\Gal(K_m /K)]$ and $\theta_g \in \Lambda(\Gamma_\ac)$ to denote the image of $\theta^{\rm KP}_{g, m}$ and $\theta^{\rm KP}_g$ under the morphisms induced by the natural projection of Galois groups. 
\end{defn}

\subsubsection{} \label{subsubsec_theta_els}
We now define the theta elements for $n$-admissible forms, and compare them with those of Definition \ref{defn_KP_theta_classes}. If $n \in \Z^{\geq 1}$, we set $\Lambda_n(\Gamma_\ac) := \Lambda(\Gamma_\ac) \otimes_\mathcal{O} \mathcal{O}_n$.
Let $\mathcal{D} = (\Delta, g)$ denote an $n$-admissible form such that $\omega(\Delta)$ is \emph{odd}. If $a\in \widehat K^\times$, recall the CM points $x_m(a)$ from Definition \ref{def_cm-points} attached to $a$, as well as Remark \ref{remark_fields_of_definition_of_CM_points}. For $m\geq n$ we set 
\begin{equation}
\Theta_m(\mathcal{D}) = \lambda_p(g)^{-m} \sum_{\sigma \in \Gal(H_m/ K)} g(x_m(a_\sigma))\cdot\sigma \in \mathcal{O}_n[\Gal(H_m/K)], \label{eqn_Big_Theta_elmt}
\end{equation}
where $a_\sigma$ denotes a lift of an element corresponding to $\sigma$ under the reciprocity map. 
By Proposition \ref{prop_distribution_relation}, we see that  $\Theta_{m+1}(\mathcal{D}) \mapsto \Theta_{m}(\mathcal{D})$ under the projection map $\mathcal{O}_n[\Gal(H_{m+1}/K)] \to \mathcal{O}_n[\Gal(H_{m}/K)]$.
\begin{defn}\label{defn_theta_elements}
    For $m\geq n$, set $\theta_m(\mathcal{D}) \in \mathcal{O}_n[\Gal(K_m/K)]$ to be the image of $\Theta_m(\mathcal{D})$ under the map induced by projection, and set $\theta(\mathcal{D}) := \varprojlim \theta_m(\mathcal{D}) \in \Lambda_n(\Gamma_\ac)$.
\end{defn}

The next proposition compares this theta element attached to the $n$-admissible form $(N^-, \widetilde f_n)$ of Example \ref{example_n-admissible_form} to the one constructed in \cite{kohen_pacetti_anticyc_l-function}, attached to our fixed form $f$.

\begin{proposition}\label{prop_comparing_thetas}
    Suppose $\mathcal{D}=\mathcal{D}_n:=(N^-, \widetilde f_n)$. Then there exists a unit $u \in \Z_p^\times$ such that $$\theta_{f, m} \equiv u\cdot \theta_m(\mathcal{D}_n) \mod \varpi^n$$
\end{proposition}

\begin{proof}
We choose $f_B$ as in Example~\ref{example_n-admissible_form} and in \S\ref{S:padicL}, taking
$\Delta=N^-$, $g=f$, $\psi=\epsilon^{-1}$ and
$g_{B_\Delta}=f_B$. We choose this common vector primitively over $\cO$ and take $\Omega_f$ with respect to this normalization.

    Recall the notation $\widetilde{G}_m = \Gal(\widetilde{H}_m/K)$, and let $G_m = \Gal(H_m/ K)$, and $F_m = \Gal(\widetilde{H}_m/H_m)$. Then we have 
    \begin{align*}
        \theta^{\rm KP}_{f,m} &= a_p(f)^{-m}\sum_{\tau \in G_m} \sum_{\sigma \in F_m} \epsilon^{-1}(\widetilde \tau \sigma)f_{B}(x_m(a_{\widetilde \tau \sigma})) \cdot \widetilde{\tau}\sigma \\
        &= a_p(f)^{-m}\sum_{\sigma \in F_m} \sigma \cdot \sum_{\tau \in G_m} \epsilon^{-1}_{\AA^\times}({\rm N}_{K/\Q} ( a_{\widetilde \tau \sigma}))f_B(x_m(a_{\widetilde \tau \sigma})) \cdot \widetilde \tau \\
        &= a_p(f)^{-m} \sum_{\sigma \in F_m} \sigma \cdot \sum_{\tau \in G_m} \Tw_{\epsilon}(f_B)\,(x_m(a_{{\tau}})) \cdot \widetilde \tau \\
        &\equiv \lambda_p(\widetilde f_n)^{-m} \sum_{\sigma \in F_m} \sigma \cdot \sum_{\tau \in G_m}\widetilde{f}_n(x_m(a_{\tau}) )\cdot \widetilde \tau \mod \varpi^n \\
        & \mapsto |F_m| \cdot \Theta_m(\mathcal{D}_n). 
    \end{align*}
    where $\widetilde \tau \in \widetilde{G}_m$ denotes a lift of $\tau\in G_m$. Here, the second equality follows from our normalization of the reciprocity map, together with the fact that we may choose $\gamma$ in the definition of $\varsigma$ (\textit{cf.} Lemma \ref{lemma_optimal_embeddings}) such that $\det \gamma \equiv 1 \mod p$. The third equality is by definition of the twisting map $\Tw_{\epsilon}$ (see Remark \ref{rk:Hecke}), together with the first sentence in Remark \ref{remark_fields_of_definition_of_CM_points}. The congruence in the fourth line is immediate by the definition of $\widetilde{f}_n$, and the map in the final line is the one induced by the projection $\widetilde G_m \to G_m$. 
    
  %  Since $F_\psi/\Q$ has degree prime to $p$, we have that $|F_m|=:u \in \Z_p^\times$, and the claim now follows after projecting to $\mathcal{O}_n[\Gal(K_m/K)]$.
Projecting to $\cO_n[\Gal(K_m/K)]$ gives
\[
\theta_{f,m}\equiv |F_m|\cdot\theta_m(\cD_n)\mod\varpi^n.
\]
Since $F_\psi/\Q$ has degree prime to $p$, we have $|F_m|\in\Zp^\times$. Moreover $H_m/H_1$ has $p$-power degree for $m\ge1$, so $F_\psi\bigcap H_m=F_\psi\bigcap H_1$. Thus, $|F_m|$ is independent of $m$, which concludes the proof.
\end{proof}

\subsection{Explicit reciprocity laws}
We record the `first and second explicit reciprocity laws', which will be used in the following section to prove the main results of this article. In what follows, we abuse notation in identifying elements of $H^1(K_{\ac, \ell}, T_n)$ with their images in $\Lambda_n(\Gamma_\ac)$ under the fixed isomorphism \eqref{eqn_isomorphism_H^1=Lambda_n}.

\begin{theorem}[First explicit reciprocity law]\label{theorem_first_erl}
    Suppose $\omega(\Delta)$ is even and there exists an $n$-admissible prime $\ell \mid \Delta$. Let $\mathcal{D} = (\Delta, g)$ be an $n$-admissible form. Then we have the equality 
    $$\partial_\ell(\kappa_{\mathcal{D}}) = \theta(\mathcal{D}^{(\ell)}) \in \Lambda_n(\Gamma_\ac).$$
\end{theorem}
\begin{proof}
    This follows as in the proof of \cite[Theorem 5.1]{chidahsiehanticyclomainconjformodformscomposito} (\textit{cf.} Lemma \ref{lemma_comparing_Heegner_classes} and Remark \ref{remark_Hida_theory}).
\end{proof}

\begin{theorem}[Second explicit reciprocity law]\label{theorem_second_erl}
   Suppose $\omega(\Delta)$ is even. Let $\ell\nmid \Delta$ be an $n$-admissible prime, and $\mathcal{D} = (\Delta, g)$ an $n$-admissible form. Then there exists an $n$-admissible form $\mathcal{D}'= (\Delta\ell, g')$ satisfying the equality
    $$\res_{\ell}(\kappa_{\mathcal{D}}) = \theta(\mathcal{D}') \in \Lambda_n(\Gamma_\ac). $$
\end{theorem}
\begin{proof}
    This follows as in the proof of \cite[Theorem 5.5]{chidahsiehanticyclomainconjformodformscomposito}, since the geometry at primes not equal to $p$ is unaffected by non-semistable reduction at $p$ of the relevant Shimura curve. The version of Ihara's lemma (Corollary 5.4 of \textit{op. cit.}) used in the proof of \textit{loc. cit.} also holds in our case by Remark \ref{remark_Hida_theory} and \eqref{item_ord}.
\end{proof}

\section{Bipartite Euler systems and specialization}
\subsection{The bipartite Euler system formalism} In this subsection, we recall and slightly extend some of the formalism developed in \cite[\S 3]{howard06}, following the exposition in \cite[\S 6]{BBL2}.

\subsubsection{}\label{subsubsec_prelims_for_main} 
Let $\mathfrak L_n$ denote the set of $n$-admissible primes, and let $\mathfrak F$ denote the Selmer structure on $T_{\widetilde{f}}$ given by the local conditions in \S \ref{subsection_Selmer_groups} (\textit{cf.} \cite[Definition 2.1.1]{howard06}).

\begin{proposition} \label{prop_Howard_hypotheses}
    Assume \eqref{item_Heeg}, \eqref{item_ord}, \eqref{item_CR} and \eqref{item_irr}. Then the triple $(T_n, \mathfrak F, \mathfrak L_n)$ satisfies the conditions of \cite[Lemma 3.3.4]{howard06}. Namely, 
    \begin{itemize}
        \item[i)] $T_n$ is residually absolutely irreducible.
        \item[ii)] For any $\ell \in \mathfrak L_n$, $\Frob_\ell^2$ acts on $T_n$ with eigenvalues $\ell^2$ and $1$.
        \item[iii)] There is a symmetric perfect $\mathcal{O}_n$-bilinear pairing $(\, ,\, ):T_n \times T_n \lra \mathcal{O}_n(1)$ satisfying $(s^\sigma, t^{\tau \sigma \tau }) = (s,t)^\sigma $ for all $\sigma \in G_K$, where $\tau$ denotes the complex conjugation. 
        \item[iv)] $\mathfrak F$ is a cartesian Selmer structure in the sense of \cite[Definition 2.2.2]{howard06} and is self-dual in the sense of \S 2.6 of $\textit{op. cit}$ for the pairing above. 
        \item[v)] For any non-zero $\kappa \in H^1(K, T_1)$, there are infinitely many primes $\ell \in \mathfrak L_n$ such that $\res_\ell(\kappa) \neq 0$.
    \end{itemize}
\end{proposition}

\begin{proof}
    (i) is Lemma \ref{lemma_Howard_hyp_a}, (ii) follows from \eqref{eqn_decomp_of_T_at_ell}. For (iii), let $\langle\,,\,\rangle$ be the alternating self-dual pairing on $T_n$ and define $(s,t)=\langle s,\tau t\rangle$. (v) is Theorem \ref{theorem_Howard_hyp_e}. For (iv), self-duality follows by Proposition \ref{prop_orthogonal_at_p}, Lemmas \ref{lemma_cohomology_decomp_at_admissible_ell} and \ref{lemma_orthogonal_complements_at_ell_dividing_N}, and \cite[Remark 2.1.2]{howard06}. The cartesian condition follows from \cite[Lemma 3.7.1]{Mazur_Rubin} and Lemma \ref{lemma_vanishing_cohomology_at_p}. 
\end{proof}

\subsubsection{} Let $\mathfrak{N}_n$ denote the set of square-free products of primes in $\mathfrak{L}_n$, together with the empty product $1$. For $\mathfrak{n} \in \mathfrak{N}_n$, we set 
\begin{align*}
    &\Delta_\mathfrak{n} := N^- \cdot \mathfrak{n}, \\
    & \mathcal{D}_{\mathfrak{n}} := (\Delta_{\mathfrak{n}}, f_{\mathfrak{n}}),
\end{align*}
so that $B_{\Delta_\mathfrak{n}}$ is definite or indefinite, depending on the parity of $\omega(\Delta_\mathfrak{n})$, and $\mathcal{D}_{\mathfrak{n}}$ is an $n$-admissible form. 
Consider the decomposition 
$$\mathfrak{N}_n = \mathfrak{N}_n^{\rm def} \bigsqcup \mathfrak{N}_n^{\rm indef},$$ where 
\begin{align*}
    \mathfrak{N}_n^{\rm def} &:= \{\mathfrak n \in \mathfrak{N}_n : \omega(\Delta_ \mathfrak{n}) \text{ is odd}\}, \\
    \mathfrak{N}_n^{\rm indef}& := \{\mathfrak n \in \mathfrak{N}_n : \omega(\Delta_ \mathfrak{n}) \text{ is even}\}.
\end{align*} 

\begin{lemma}\label{lem:datum}
    With notation as above, $\mathcal{D}_{\mathfrak{n}}$ exists.
\end{lemma}

\begin{proof}
    This follows by induction on the prime factors of $\mathfrak{n}$. For $\mathfrak{n}=1$ we take $\mathcal{D}_1$ to be as in Example \ref{example_n-admissible_form}. Now suppose we are given an $n$-admissible form $\mathcal{D}_\mathfrak{n}$. If $\mathfrak{n} \in \mathfrak{N}^{\rm indef}_n$, then Theorem \ref{theorem_second_erl} gives an $n$-admissible form for $\mathfrak{n}\ell \in \mathfrak{N}^{\rm def}_n$. On the other hand, if $\mathfrak{n} \in \mathfrak{N}^{\rm def}_n$, then we can construct an $n$-admissible form for $\mathfrak{n}\ell \in \mathfrak{N}^{\rm indef}_n$ as follows. First let $\T^{[\ell]}_{\Delta_{\mathfrak{n}}}$ denote the Hecke algebra given by replacing $T_\ell$ by $U_\ell$ and consider the morphism $f_{\mathfrak{n}}^{[\ell]}: \T^{[\ell]}_{\Delta_{\mathfrak{n}}} \lra \mathcal{O}_n$  given by $f_\mathfrak{n}$ away from $U_\ell$, and $U_\ell \mapsto \varepsilon_\ell$. By \cite[Theorem 4.3]{chidahsiehanticyclomainconjformodformscomposito}, the map $f^{[\ell]}_\mathfrak{n}$ factors through the map $\T^{[\ell]}_{\Delta_{\mathfrak{n}}} \twoheadrightarrow \T_{\Delta_{\mathfrak{n}}\ell}$ given in \S3.4 of \textit{op. cit.} and hence induces an $n$-admissible form $\mathcal{D}_{\mathfrak{n}\ell}$.
\end{proof}

Fix a discriminant $\Delta$ and write 
\[
X_s=X(\su_s),\quad J_s=J(\su_s),\quad \pi_{n,r}:X_n\to X_r
\]
for $1\le r\le n$. Let $\TT_{\Delta,s}$ be the Hecke algebra given in \S\ref{S:hecke} for the level $\su_s$, and write $e_\ord=\displaystyle\lim_{n\rightarrow\infty} U_p^{n!}$. 

In the definite case, for $R=\cO$ or a quotient of $\cO$, define
\[
\cS_s^\epsilon(R)=e_\ord\left\{h\in M_2^\Delta(\su_s,R):\langle a\rangle h=\epsilon(a)h, a\in\Zp^\times\right\}.
\]

In the indefinite case, let
$g_s:\mathbb T_{\Delta,s}\to\mathcal O_s$
be the eigencharacter of the chosen $s$-admissible form. Define $ \mathfrak m_s$ as the maximal ideal of $\TT_{\Delta,s}$ given by $\ker(\TT_{\Delta,s}\stackrel{g_s}\to\cO_s\to \cO_1)$.
Put
\[
\cH_s=e_\ord(T_p(J_s)\otimes\cO)_{\fM_s},\quad P_{\epsilon,s}=(\langle a\rangle-\epsilon(a):a\in\Zp^\times)\subset (\TT_{\Delta,s})_{\fM_s},\quad\cM_s=\cH_s/P_{\epsilon,s}\cH_s.
\]
\begin{lemma}\label{lem:control}
The following properties hold.
\begin{itemize}
    \item[i)]In the definite case, we have an isomorphism $\pi_{n,r}^*:\cS_r^\epsilon(R)\stackrel{\sim}\longrightarrow\cS_n^\epsilon(R)$.
    \item[ii)] In the indefinite case, we have an isomorphism $(\pi_{n,r})_*:\cM_n\stackrel{\sim}\longrightarrow\cM_r$.
    \item[iii)]The isomorphisms above are Hecke-equivariant, commute with coefficient reduction, and are transitive.  
\end{itemize}
    \end{lemma}
    \begin{proof}
        Undoing the twist in Remark~\ref{rk:Hecke} reduces the comparison to the ordinary parts of the spaces at the levels $\su_s^1$ defined in Remark~\ref{remark_Comparing_level_subgroups}, with weight two and nebentype $\epsilon^2$. By Remark~\ref{remark_Hida_theory}, we can apply Hida's ordinary control theorem to these spaces; see \cite[Corollary~8.2 and Proposition~8.3]{hid88} and \cite[proofs of Proposition 4.2 and Corollary 5.4]{chidahsiehanticyclomainconjformodformscomposito}. The twisting factors are $p$-adic units, so these identifications are integral.
    \end{proof}

Let us write $\cD_\fn^{(n)}=\left(\Delta_\fn,f_\fn^{(n)}\right)$ for the admissible datum given in Lemma~\ref{lem:datum}. For $\fn\ne1$, set $h(\fn)=\max\{s\ge1:\fn\in\fN_s\}$.
The integer $h(\fn)$ is finite: the Weil
bound gives
\[
 |\epsilon^{-1}(\ell)a_\ell(f)|
 \leq 2\sqrt{\ell}<\ell+1
 \qquad(\ell\mid\fn),
\]
so $(\ell+1)^2-\epsilon^{-2}(\ell)a_\ell(f)^2\neq0$,
and the assertion follows from Definition~\ref{def_n-admissible_primes}(iv).

For every $\fn\ne1$, we can choose the admissible datum  for $n=h(\fn)$, and define the corresponding datum for $r\le n$ by reduction through Lemma~\ref{lem:control}. For $\fn=1$, we use Example~\ref{example_n-admissible_form}. In the definite case, this means
$f_{\fn}^{(n)}\mod\varpi^r=\pi_{n,r}^*f_\fn^{(r)}$ for $1\le r\le n$. In the indefinite case, the corresponding eigencharacters are compatible under reduction.

Suppose $\fn$ is indefinite. For each $r$ such that $\fn\in\fN_r$, set
\[
Q_{\fn,r}:=(T_p(J(\su_r))\otimes \cO)/I_{f_\fn^{(r)}}\cong\cM_r/I_{f_\fn^{(r)}}.
\]
This quotient is identified with $T_r$ in \S\ref{s:heegnerclasses}.
By coefficient compatibility and the identification of Hecke algebras of Lemma~\ref{lem:control}(ii), for $1\le r\le n$, $I_{f_\fn^{(r)}}=I_{f_\fn^{(n)}}+(\varpi^r)$, and we have $$\rho_{n,r}:Q_{\fn,n}\to Q_{\fn,r},\quad Q_{\fn,n}/\varpi^r Q_{\fn,n}\stackrel{\sim}{\longrightarrow} Q_{\fn,r}.$$
Let $\re_{n,r}:T_n\to T_r$ be the reduction map. We can choose $\cO_r[G_K]$-isomorphisms $\phi_r:Q_{\fn,r}\stackrel{\sim}{\longrightarrow} T_r$ such that $\phi_r\circ\rho_{n,r}=\re_{n,r}\circ\phi_n$. Indeed, Lemma~\ref{lemma_Howard_hyp_a} and Schur's lemma give $\mathrm{Aut}_{\cO_r[G_K]}(T_r)=\cO_r^\times$. After choosing $\phi_r$, we can choose any  $Q_{\fn,r+1}\stackrel{\sim}\lra T_{r+1}$ and rescale
it by a unit so that its reduction agrees with $\phi_r$.

Let $x_m^{(r)}(a)$ and $D_m^{(r)}$ be the CM point and divisor defined for the level $\su_r$ in Definition~\ref{def_cm-points} and \eqref{eq:CM}, respectively. Suppose that the same $q_0$ and CM representatives are used for each $\su_r$ as $r$ varies. For $r\le n\le m$, we have
\[
\pi_{n,r}(x_m^{(n)}(a))=x_m^{(r)}(a),\quad (\pi_{n,r})_*D_m^{(n)}=D_m^{(r)}.
\]
By functoriality of the Kummer maps and evaluation at CM points, we have
\begin{align*}
    \kappa_{\cD_\fn^{(n)}}\mod \varpi^r&=\kappa_{\cD_\fn^{(r)}},\\
\theta(\cD_\fn^{(n)})\mod\varpi^r&=\theta(\cD_\fn^{(r)}),
\end{align*}
whenever these elements are defined.

\begin{defn}[Bipartite Euler system] \label{defn_bipartite_Euler_system}
    A bipartite Euler system for $(T_{\widetilde{ f}}, \mathfrak{F}, (\mathfrak{L}_n)_n)$ is the data of two families of elements 
    \begin{align*}
        &\{\kappa(\mathfrak{n},n) \in \Sel_{\Delta_\mathfrak{n}}(K_\ac, T_n) : \mathfrak{n} \in \mathfrak{N}^{\rm indef}_n\}, &\{\theta(\mathfrak{n}, n) \in \Lambda_n(\Gamma_\ac) : \mathfrak{n} \in \mathfrak{N}^{\rm def}_n \},  
    \end{align*}
    for all positive integers $n$, such that the following conditions hold: 
    \begin{itemize}
        \item[i)]  The elements $\kappa(\mathfrak{n},n)$, $\theta(\mathfrak{n}, n)$ are compatible under the natural maps as $n$ varies, i.e., 
        \[
        \kappa(\fn,n)\mod \varpi^r=\kappa(\fn,r),\qquad \theta(\fn,n)\mod \varpi^r=\theta(\fn,r)
        \]
        for $r\le n$ and $\fn\in\fN_n$.
        \item[ii)] For any $\mathfrak{n} \in \mathfrak{N}_n^{\rm indef}$, with $\ell \mid \mathfrak{n}$ an $n$-admissible prime, there is an isomorphism of $\Lambda(\Gamma_\ac)$-modules 
        \[H^1_{\rm ord}(K_{\ac, \ell}, T_n) \cong \Lambda_n(\Gamma_\ac)\]
        which carries $\res_\ell(\kappa(\mathfrak{n}, n))$ to $\theta(\mathfrak{n}/\ell, n)$. 
        \item[iii)] For any $\mathfrak{n}\ell \in \mathfrak{N}_n^{\rm def}$, there is an isomorphism of $\Lambda(\Gamma_\ac)$-modules 
        \[H^1_{\rm fin}(K_{\ac, \ell}, T_n) \cong \Lambda_n(\Gamma_\ac)\]
        which carries $\res_\ell(\kappa(\mathfrak{n}, n))$ to $\theta(\mathfrak{n}\ell, n)$. 
    \end{itemize}
\end{defn}
 
\begin{lemma} \label{lemma_Bipartite_Euler_System_Exists}
    The sets 
    \begin{align*}
        &\{\kappa_{\mathcal{D}_\mathfrak{n}} \in \Sel_{\Delta_\mathfrak{n}}(K_\ac, T_n) : \mathfrak{n} \in \mathfrak{N}^{\rm indef}_n\}, &\{\theta(\mathcal{D}_\mathfrak{n}) \in \Lambda_n(\Gamma_\ac) : \mathfrak{n} \in \mathfrak{N}^{\rm def}_n \},  
    \end{align*}
    form a bipartite Euler system for $(T_{\widetilde{f}}, \mathfrak{F}, (\mathfrak{L}_n)_n)$, with $\kappa(\mathfrak{n},n) = \kappa_{\mathcal{D}_{\mathfrak{n}}}$, $\theta(\mathfrak{n}, n) = \theta(\mathcal{D}_\mathfrak{n})$. 
    %Moreover, we may shrink each $\mathfrak{L}_n$ as necessary so that the above bipartite Euler system is free in the sense of \cite[Definition 2.3.6]{howard06}.
\end{lemma}

\begin{proof}
    The compatibility of  $\kappa_{\cD_\fn}$ and $\theta(\cD_\fn)$ discussed above gives the property asserted in Definition~\ref{defn_bipartite_Euler_system}(i). Theorems \ref{theorem_first_erl} and \ref{theorem_second_erl} give (ii) and (iii), with the asserted isomorphisms normalized to match our chosen forms and
Tate-module identifications. %The second claim follows from the same argument as in \cite[Lemma 3.3.6]{howard06}, with the relevant hypotheses being provided by Proposition \ref{prop_Howard_hypotheses}.
\end{proof}

%We henceforth always assume that the bipartite Euler system of the previous lemma is free.

\subsubsection{} 
Let us set 
\begin{align*}
    &\mathbf{T} := \varprojlim_m \Ind^{K_m}_{K} T_{\widetilde{f}}, &\mathbf{A} := \varinjlim_m \Ind^{K_m}_{K} A_{\widetilde{f}},
\end{align*}
so that $\mathbf{T} \cong T_{\widetilde{f}} \otimes_{\mathcal{O}} \Lambda(\Gamma_\ac)$, where $G_{K}$ acts on both factors. Let $\mathfrak{P} \subset \Lambda(\Gamma_\ac)$ be a height-one prime ideal with $\mathfrak{P} \neq \varpi \cdot \Lambda(\Gamma_\ac)$, and write $\mathcal{O}_{\mathfrak{P}}$ for the integral closure of $\Lambda(\Gamma_\ac)/\mathfrak{P}$. We set 
\begin{align*}
    &T_{\mathfrak{P}}: = \mathbf{T} \otimes \mathcal{O}_{\mathfrak{P}}, &A_{\mathfrak{P}}:= T_{\mathfrak{P}} \otimes \Q_p/ T_{\mathfrak{P}}.
\end{align*}
Moreover, if $\mathfrak{P}^\iota$ denotes the height one prime ideal obtained from $\mathfrak{P}$ via the canonical involution on $\Gamma_\ac$, then as in \cite{howard04}, since $T_{\widetilde{f}}$ is self-dual, we have a natural pairing $T_{\widetilde{f}} \times A_{\widetilde{f}} \lra \mu_{p^{\infty}}$, which induces a perfect $G_K$-equivariant pairing 
\begin{equation}
    T_{\mathfrak{P}^\iota} \times A_{\mathfrak{P}} \lra \mu_{p^\infty} \label{eqn_perfect_pairing}.
\end{equation}
The Selmer structure $\mathfrak{F}$ (\textit{cf.} \S \ref{subsubsec_prelims_for_main}) induces a Selmer structure on $\mathbf{T}$, $\mathbf{A}$. On $T_{\mathfrak P}$ and $A_{\mathfrak P}$, we use the
inverse images and images, respectively, of the rational
local conditions, ordinary at $pN^-$ and unramified
elsewhere. The conditions on finite subquotients are
obtained by propagation. We continue to denote the induced Selmer groups by $\Sel_{\Delta}$. In particular, if $v\mid p$ is a prime of $K$, we can consider the induced local conditions 
\begin{align*}
    &H^1_{\ord}(K_v, \mathbf{T}/ \mathfrak{P}\mathbf{T}), &H^1_{\ord}(K_v, \mathbf{A}[\mathfrak{P}]), && H^1_{\ord}(K_v, {T}_{\mathfrak{P}}), &&H^1_{\ord}(K_v, {A}_{\mathfrak{P}}).
\end{align*}

\begin{lemma}\label{lem:ord-transition}
    Assume \eqref{item_Heeg}, \eqref{item_ord} hold. Let $\mathfrak{P} \neq \varpi \cdot \Lambda(\Gamma_\ac)$ be a height-one prime and let $v$ be a place of $K$ dividing $pN^-$. Then the natural maps 
    \begin{align}
        &H^1_{\ord}(K_v, \mathbf{T}/\mathfrak{P}\mathbf{T}) \lra H^1_{\ord}(K_v, T_{\mathfrak{P}}), \label{eqn_local_map_1} \\
        &H^1_{\ord}(K_v, A_\mathfrak{P}) \lra H^1_{\ord}(K_v,\mathbf{A}[\mathfrak{P}]) \label{eqn_local_map_2},
    \end{align}
    where the second map is induced by \eqref{eqn_perfect_pairing}, have finite kernels and cokernels whose orders are bounded by constants depending only on $[\mathcal{O}_{\mathfrak{P}}: \Lambda(\Gamma_\ac)/\mathfrak{P}]$.
\end{lemma}
\begin{proof}
    The same proof of \cite[Lemma 2.2.7]{howard04} applies verbatim, with some simplifications in our case, using Lemma \ref{lemma_vanishing_cohomology_at_p}. 
\end{proof}

\begin{proposition} \label{prop_control_theorems}
    The natural maps 
    \begin{align*}
    \Sel_{N^-}(K, \mathbf{T})/\mathfrak{P} \,\Sel_{N^-}(K, \mathbf{T}) &\lra \Sel_{N^-}(K, T_{\mathfrak{P}}), \\
        \Sel_{N^-}(K, A_{\mathfrak{P}})&\lra \Sel_{N^-}(K,\mathbf{A} )[\mathfrak{P}],
    \end{align*}
    satisfy the following conditions: 
    \begin{itemize}
        \item[i)] The first map is injective.
        \item[ii)] Outside a finite set of height one primes, both maps have finite kernel and cokernel bounded by a constant depending only on $[\mathcal{O}_{\mathfrak{P}}:\Lambda(\Gamma_{\ac})/\mathfrak{P}]$.
    \end{itemize}
\end{proposition}
\begin{proof}
    The same proof as in \cite[Proposition 3.3.1]{howard06} applies here, where we use the previous lemma in place of \cite[Lemma 2.2.7]{howard04} in the proof of \textit{loc. cit.}
\end{proof}

\begin{lemma}\label{lem:non-vanishing}
    The natural map 
    \begin{equation*}
        \Sel_{N^-}(K, T_{\mathfrak{P}})/ \varpi^n \Sel_{N^-}(K, T_{\mathfrak{P}}) \lra \Sel_{N^-}(K, T_{\mathfrak{P}}/\varpi^n T_{\mathfrak{P}})
    \end{equation*}
    is injective, where the Selmer structure on $T_{\mathfrak{P}}/\varpi^n T_{\mathfrak{P}}$ is obtained from that of $T_{\mathfrak{P}}$ by propagation.
\end{lemma}
\begin{proof}
    The same proof as in \cite[Lemma 3.3.2]{howard06} applies, which in turn is a special case of \cite[Lemma 3.7.1]{Mazur_Rubin}. 
\end{proof}

\subsubsection{} To ease notation, we put 
\begin{align*}
    &\mathcal{O}_{\mathfrak{P},n} := \mathcal{O}_{\mathfrak{P}}/\varpi^n \mathcal{O}_{\mathfrak{P}}, \quad T_{\mathfrak{P},n} := T_{\mathfrak{P}}/\varpi^n T_{\mathfrak{P}},\quad e_\fP=\ord_{\fM_\fP}(\varpi),
\end{align*}
where $\fM_\fP$ is the maximal ideal of $\cO_\fP$.

A choice of uniformizer for $\mathcal{O}_{\mathfrak{P}}$ determines an isomorphism $A_{\mathfrak{P}}[\varpi^n] \cong T_{\mathfrak{P}, n}$ which in turn gives us the identification
\begin{equation}
    \Sel_{N^-}(K, T_{\mathfrak{P}, n}) \cong \Sel_{N^-}(K, A_{\mathfrak{P}}[\varpi^n]) \cong \Sel_{N^-}(K, A_{\mathfrak{P}})[\varpi^n], \label{eqn_sel_torsion}
\end{equation}
by \cite[Lemmas 3.7.1 \& 3.5.4]{Mazur_Rubin}. If $M$ is an $\mathcal{O}_{\mathfrak{P}}$-module, and $m \in M$, then the index of divisibility of $m$, denoted $\Ind(m, M )$, is the largest integer $r \leq \infty $ such that $m \in \mathfrak{m}_{\mathfrak{P}}^rM$. For positive integers $k\leq n$, set  
\begin{equation*}
    \delta_{\mathfrak{P}}(k,n) := \min \{\Ind(\theta(\mathcal{D}_{\mathfrak{n}}), \mathcal{O}_{\mathfrak{P},k}): \mathfrak{n} \in \mathfrak{N}^{\rm def}_n \}. 
\end{equation*}
Note that $\delta_{\mathfrak{P}}(k,n)$ increases monotonically in $n$, and define 
\begin{equation*}
    \delta_{\mathfrak{P}}(k) := \lim_{n \to \infty} \delta_{\mathfrak{P}}(k, n).
\end{equation*}

\subsubsection{}
Suppose momentarily that $\omega(N^-)$ is odd. We may consider the theta elements $\theta(\mathcal{D}^{(n)}_1) \in \Lambda_n(\Gamma_\ac)$ attached to the $n$-admissible forms $\mathcal{D}^{(n)}_1 := (N^-, \widetilde f_n)$ of Example \ref{example_n-admissible_form}. These are compatible under the natural maps by Lemma \ref{lemma_Bipartite_Euler_System_Exists}, and hence we may define
$$\theta_\infty := \varprojlim_n \theta(\mathcal{D}^{(n)}_1) \in \Lambda(\Gamma_\ac).$$
If instead $\omega(N^-)$ is even, then we may consider the Heegner classes 
$$\kappa_{\mathcal{D}^{(n)}_1} \in \Sel_{N^-}(K_\ac, T_n),$$
attached to the $n$-admissible forms $\mathcal{D}^{(n)}_1$. These classes are again compatible under the natural morphisms by Lemma \ref{lemma_Bipartite_Euler_System_Exists} and hence we may define the class
\begin{equation}\label{eq:kappa_infty}
\kappa_\infty := \varprojlim_n \kappa_{\mathcal{D}^{(n)}_1} \in \Sel_{N^-}(K, \mathbf T).    
\end{equation}

\begin{proposition} \label{prop_lengths_of_Selmer}
    \item[i)] Suppose $\omega(N^-)$ is odd, so we are in the definite case, and assume that there exists a positive integer $k$ such that $\theta_\infty \in \Lambda(\Gamma_\ac)$ has non-trivial image in $\mathcal{O}_{\mathfrak{P}, k}$. Then we have the equality
    \begin{equation*}
        {\rm length}_{\mathcal{O}_{\mathfrak{P}}}\left( \Sel_{N^-}(K, A_{\mathfrak{P}})\right) + 2\delta_{\mathfrak{P}}(k) = 2 \cdot {\rm length}_{\mathcal{O}_{\mathfrak{P}}}\left( \mathcal{O}_{\mathfrak{P}}/ \theta_\infty \right).
    \end{equation*}
    \item[ii)] Suppose $\omega(N^-)$ is even, so we are in the indefinite case, and assume that there exists a positive integer $k$ such that $\kappa_{\mathfrak{P}} \not \equiv 0 \mod \varpi^k$, where $\kappa_{\mathfrak{P}}$ denotes the image of $\kappa_\infty \in \Sel_{N^-}(K, \mathbf{T})$ in $\Sel_{N^-}(K, T_{\mathfrak{P}})$. Then the following hold: 
    \begin{itemize}
        \item[a)] The $\mathcal{O}_{\mathfrak{P}}$-module $\Sel_{N^-}(K, T_{\mathfrak{P}})$ is free of rank one.
        \item[b)] The $\mathcal{O}_{\mathfrak{P}}$-module $\Sel_{N^-}(K, A_{\mathfrak{P}})^\vee$ is of rank one.
        \item[c)] We have the equality
        \begin{equation*}
            {\rm length}_{\mathcal{O}_{\mathfrak{P}}}\left( \Sel_{N^-}(K, A_{\mathfrak{P}})^\vee_{\rm tor} \right) + 2\delta_{\mathfrak{P}}(k) = 2 \cdot {\rm length}_{\mathcal{O}_{\mathfrak{P}}}\left( \Sel_{N^-}(K, T_{\mathfrak{P}})/ \kappa_{\mathfrak{P}}\right),
        \end{equation*}
        where the subscript ${\rm tor}$ means the maximal torsion subgroup. 
    \end{itemize}
\end{proposition}

\begin{proof}
    %We proceed as in the discussion at the end of \cite[\S 3.3]{howard06}. When $\omega(N^-)$ is odd, the desired equality follows directly by \cite[Theorem 2.5.1]{howard06}, which applies thanks to Proposition \ref{prop_Howard_hypotheses}.
Fix $n\geq2k$. Note that
\[
 T_{\fP}/\fM_{\fP}\cong T_1\otimes_{\cO_1}
       (\cO_{\fP}/\fM_{\fP})
\]
is absolutely irreducible.
At primes above $N^-$ and in $\fL_n$, which split completely in $K_{\mathrm{ac}}/K$, the representation $T_{\mathfrak P}$ is obtained from $T_{\widetilde f}$ by scalar extension. Hence Proposition~\ref{prop_Howard_hypotheses}(i), (ii) and (v) apply after residue-field extension. The pairing in (iii) also extends, since complex conjugation acts on $\Gamma_{\mathrm{ac}}$ by inversion.

For $v\mid p$, the short exact sequence \eqref{eq:filtration} induces
\[
 0\lra F_p^+T_{\fP}  \lra T_{\fP}  \lra F_p^-T_{\fP}\lra0.
\]
Write $F_p^\pm A_{\fP}$ for the corresponding submodule and quotient of $A_{\fP}$.
The reductions of $F_p^\pm T_{\fP}$ are scalar
extensions of $F_p^\pm T_1$, so the argument of
the proof of Lemma~\ref{lemma_vanishing_cohomology_at_p} gives
\[
 H^0(K_v,F_p^\pm A_{\fP})
 =H^2(K_v,F_p^\pm T_{\fP})=0.
\]
These vanishings identify the ordinary and propagated local conditions at $p$. Away from $p$, the comparisons follow as in Proposition~\ref{prop_Howard_hypotheses} and \cite[Proposition~1.9(2)]{chidahsiehanticyclomainconjformodformscomposito}, with the corrected hypotheses of \cite{CH-erratum}.
The propagated conditions are cartesian by \cite[Lemma~3.7.1]{Mazur_Rubin} and self-dual by local Tate duality. Thus, \cite[Lemma~3.3.4]{howard06} applies to $(T_{\fP,r},\fF,\fL_n)$ for $r=k,n$.

Let $\ES(n)$ be obtained from Lemma~\ref{lemma_Bipartite_Euler_System_Exists} by the natural map  $\Lambda_n(\Gamma_{\ac})\lra\mathcal O_{\fP,n}$.
Define $\ES(k)$ similarly, and let
$\ES(n)_k$ be the reduction of $\ES(n)$ modulo $\varpi^k$. Compatibility gives
\[
 \ES(n)_k=\ES(k)|_{\fN_n}.
\]
The nonvanishing hypothesis, together with Lemma~\ref{lem:non-vanishing} in the indefinite case, makes $\ES(n)$ and $\ES(n)_k$ nonzero.
By the argument of \cite[Lemma~3.3.5]{howard06}, the definite and indefinite indices agree with the even and odd indices.
Since $n\geq2k$, the compatible lifts give freeness of $\mathrm{ES}(n)_k$ by \cite[Lemma~3.3.6]{howard06};
the coefficient lengths satisfy $e_{\mathfrak P}n\geq2e_{\mathfrak P}k$. We may therefore apply \cite[Theorem~2.5.1]{howard06}, and part~(i) follows as in Proposition~3.3.3 of \textit{op. cit}.

    Now suppose $\omega(N^-)$ is even. By \eqref{eqn_sel_torsion}, we have the following isomorphism
    \[
    \Sel_{N^-}(K,T_\fP)\cong\varprojlim_r \Sel_{N^-}(K,A_\fP)[\varpi^r],
    \]
    where the connecting maps are multiplication by $\varpi$.
    Let us write
    \[
    \Sel_{N^-}(K,A_\fP)\cong(\cO_\fP\otimes\Qp/\Zp)^d\oplus C,
    \]
    where $C$ is finite. Then, 
    \[
    \Sel_{N^-}(K,T_\fP)\cong\cO_\fP^d.
    \]
    Since $\kappa_\fP\in \Sel_{N^-}(K,T_\fP)$ is assumed to be non-zero, we have $d\ge1$.

    Let $\kappa_{\fP,k}$ denote the image of $\kappa_\fP$ in $\Sel_{N^-}(K,T_{\fP,k})$, which is non-zero by Lemma~\ref{lem:non-vanishing}. For $n\ge2k$, the free Euler system obtained above combined with \cite[Theorem~2.5.1]{howard06} gives
    \[
    \Sel_{N^-}(K,T_{\fP,k})\cong \cO_{\fP,k}\oplus M_k\oplus M_k,
    \]
    where 
    \[
    \length_{\cO_\fP}(M_k)+\delta_\fP(k,n)=\Ind(\kappa_{\fP,k},\Sel_{N^-}(K,T_{\fP,k}))<e_\fP k.
    \]
    In particular, $\fM_\fP^{e_\fP k-1}M_k=0$.
By \eqref{eqn_sel_torsion}, we have
\[
\Sel_{N^-}(K,T_{\fP,k})\cong\cO_{\fP,k}^d\oplus C[\varpi^k]. 
\]
Thus, we deduce that $d=1$, which proves (a) and (b).

We also have $\fM_\fP^{e_\fP k-1}C=0$ and $C\cong M_k\oplus M_k$.
As the image of $\Sel_{N^-}(K,T_\fP)/\varpi^k$ in $\Sel_{N^-}(K,A_\fP)[\varpi^k]$ under \eqref{eqn_sel_torsion} is $\Sel_{N^-}(K,A_\fP)_{\mathrm{div}}[\varpi^k]$, we have 
\[
\Ind(\kappa_{\fP,k},\Sel_{N^-}(K,T_{\fP,k}))=\length_{\cO_\fP}(\Sel_{N^-}(K,T_\fP)/\kappa_\fP).
\]
Consequently, for $n\ge 2k$,
\begin{align*}
    \length_{\cO_\fP}\left(\Sel_{N^-}(K,A_\fP)^\vee_{\mathrm{tor}}\right)+2\delta_\fP(k,n)=2\length_{\cO_\fP}(\Sel_{N^-}(K,T_\fP)/\kappa_\fP).
\end{align*}
Taking $n\rightarrow\infty$ proves (c).
\end{proof}

\section{Main conjectures}
\subsection{The definite case}

 \subsubsection{The definite main conjecture} In this subsection we will assume that $\omega(N^-)$ is odd.

\begin{proposition} \label{prop_theta_not_trivial}
    Assume \eqref{item_Heeg}, \eqref{item_ord} and $\omega(N^-)$ is odd. Then $\theta_f \neq 0$. 
\end{proposition}

\begin{proof}
    This follows immediately by combining Theorem \ref{theorem_p-adic_l-function_interpolation} and \cite[Theorem 1.4]{cornutvatsal2007}.
\end{proof} 
 
\begin{theorem}\label{thm_definite_main_conj}
    Assume \eqref{item_Heeg}, \eqref{item_ord}, \eqref{item_CR}, \eqref{item_irr} and that $\omega(N^-)$ is odd. Then 
    \begin{itemize}
        \item[i)] $\Sel_{N^-}(K_\ac, T_{\widetilde{f}}) = 0$ and $\Sel_{N^-}(K_\ac, A_{\widetilde{f}})^{\vee}$ is a torsion $\Lambda(\Gamma_\ac)$-module.
        \item[ii)] We have the divisibility $$\left.\Char_{\Lambda(\Gamma_\ac)}\left(\Sel_{N^-}(K_\ac, A_{\widetilde{f}})^\vee \right) \right| (\mathscr{L}_f).$$
        \item[iii)] Equality holds in (ii) if for every height-one prime $\mathfrak{P} \subset\Lambda(\Gamma_\ac)$, there exists an integer $k_0$ depending on $\mathfrak{P}$ such that for all $n\geq k_0$, there exists $\mathfrak{n} \in \mathfrak{N}_n^{\rm def}$ such that $\theta(\mathcal{D}_n) \not \equiv 0 \mod (\mathfrak{P},\varpi^{k_0}). $
    \end{itemize}
\end{theorem}

\begin{proof}
By Propositions \ref{prop_comparing_thetas} and    \ref{prop_theta_not_trivial}, $\theta_\infty\ne 0$. Applying the argument of \cite[\S 3.4]{howard06} to Propositions \ref{prop_control_theorems} and \ref{prop_lengths_of_Selmer} gives the assertion in part (i),
\[
\left.\Char_{\Lambda(\Gamma_\ac)}\left(\Sel_{N^-}(K_\ac, A_{\widetilde{f}})^\vee \right) \right| (\theta_\infty^2),
\]
and equality holds under the criterion stated in part (iii). 

Let $(F)=\Char_{\Lambda(\Gamma_\ac)}\left(\Sel_{N^-}(K_\ac, A_{\widetilde{f}})^\vee \right)$. We have $(F^*)=(F)$ since the complex conjugation preserves the local conditions and induces an isomorphism
\[
\Sel_{N^-}(K_\ac, A_{\widetilde{f}})^\iota\cong \Sel_{N^-}(K_\ac, A_{\widetilde{f}}),
\]
where $\iota$ denotes the inversion map on $\Gamma_\ac$. Therefore, $(F)|(\theta_\infty^*)^2$, and so $(F)^2|(\theta_\infty\theta_\infty^*)^2=(\mathscr L_f)^2$. Hence, $(F)|(\mathscr L_f)$ by unique factorization in $\Lambda(\Gamma_\ac)$, proving (ii). Under the hypothesis of (iii), this divisibility becomes an equality.
\end{proof}

\subsection{The indefinite case}
In this subsection, we assume that $\omega(N^-)$ is even and \eqref{item_irr} holds.

\subsubsection{The indefinite main conjecture} \label{section_indef_main_conj}
We will prove the following proposition in the next subsubsection.

\begin{proposition} \label{prop_Heegner_nonvanishing_mod_pi}
    Assume \eqref{item_Heeg}, \eqref{item_ord}, \eqref{item_CR}, \eqref{item_irr} and $\omega(N^-)$ is even. Then $$\kappa_\infty \not \equiv 0 \mod \varpi.$$ 
\end{proposition}

\begin{theorem} \label{them_indef_main_conjecture}
    Assume \eqref{item_Heeg}, \eqref{item_ord}, \eqref{item_CR}, \eqref{item_irr} and $\omega(N^-)$ is even. Then 
    \begin{itemize}
        \item[i)] $\Sel_{N^-}(K_\ac, T_{\widetilde f})$ is a torsion-free $\Lambda(\Gamma_\ac)$-module and $${\rm rank}_{\Lambda(\Gamma_\ac)} \left( \Sel_{N^-}(K_\ac, T_{\widetilde f}) \right) = {\rm rank}_{\Lambda(\Gamma_\ac)} \left( \Sel_{N^-}(K_\ac, A_{\widetilde f})^\vee \right)=1.$$
        \item[ii)] We have the divisibility
        $$\left.\Char_{\Lambda(\Gamma_\ac)}\left(\Sel_{N^-}(K_\ac, A_{\widetilde{f}})_{\rm tor}^\vee \right) \right| \Char_{\Lambda(\Gamma_\ac)}\left(\Sel_{N^-}(K_\ac, T_{\widetilde{f}})/ \kappa_{\infty} \right)^2 .$$
        \item[iii)] Equality holds in (ii) if for every height-one prime $\mathfrak{P} \subset\Lambda(\Gamma_\ac)$, there exists an integer $k_0$ depending on $\mathfrak{P}$ such that for all $n\geq k_0$, there exists $\mathfrak{n} \in \mathfrak{N}_n^{\rm def}$ such that $\theta(\mathcal{D}_n) \not \equiv 0 \mod (\mathfrak{P},\varpi^{k_0}). $
    \end{itemize}
\end{theorem}
\begin{proof}
    By Proposition \ref{prop_Heegner_nonvanishing_mod_pi}, this again follows by the same argument given in \cite[\S 3.4]{howard06}, where we now use Propositions \ref{prop_control_theorems} \& \ref{prop_lengths_of_Selmer} (ii).
\end{proof}

\subsubsection{Non-vanishing of the residual Heegner class} The goal of this subsubsection is to prove Proposition \ref{prop_Heegner_nonvanishing_mod_pi}. By the second explicit reciprocity law (Theorem \ref{theorem_second_erl}), it suffices to show that a theta element is non-vanishing modulo $\varpi$, which we prove following the ideas in \cite{Vatsal2003, cornutvatsal_2005, cornutvatsal2007, chidahseihanticycloLvaluescrelle}.

Let $\ell$ be a $1$-admissible prime (\textit{cf}. Theorem \ref{theorem_Howard_hyp_e}) and let 
$$B' := B_{N^-\ell},$$ 
so that $B'$ is a definite quaternion algebra and $U' := \mathscr{U}_1 \subset \widehat B'^\times $ is the compact open subgroup defined as in \S\ref{subsub_QuaternionPrelims}. By Theorem \ref{theorem_second_erl}, there exists a $1$-admissible form $\mathcal{D}' = (N^-\ell, g)$ such that 
\begin{equation*}
    \res_\ell(\kappa_{\mathcal{D}}) =\theta(\mathcal{D}') \in \Lambda_1(\Gamma_\ac), 
\end{equation*}
where $\mathcal{D} = (N^-, \widetilde{f}_1)$. By definition, $g \in \Hom(X(U'), \mathcal{O}_1)$ satisfies $T_q(g) = \lambda_q(g) g$ with $\lambda_q(g) \equiv \epsilon^{-1}(q)a_q(f) = \Tr \, {\overline \rho_{\widetilde{f}}}(\Frob_q) \mod \varpi $ for $q\nmid pN\ell$. It remains to show that $\theta(\mathcal{D}')$ is not trivial.

\begin{defn}
    Let $N: \widehat B'^\times \lra \widehat \Q^\times $ denote the reduced norm map. A function $h: X(U') \lra \mathcal{O}_1$ is Eisenstein if $h$ is invariant under the action of $\ker(N)$ induced by right-translation. 
\end{defn}

We continue to denote by $N$ the map $X(U') \lra Z:= \Q_{>0}^\times \backslash \widehat \Q^\times / N(U')$ induced by the reduced norm. Since $N(U'_p) = (\Z_p^\times)^2$ and $N(U'_q) = \Z_q^\times$ for $q\neq p$, we see that $Z$ is a group of order $2$ and we have an identification 
\begin{align*}
    Z &\xrightarrow{\sim} \{\pm 1\} \\
    z &\mapsto (\frac{z_p}{p}), \ ( z \in \widehat \Z^\times).  
\end{align*}

\begin{lemma} \label{lemma_g_not_Eis}
    Assume \eqref{item_CR_2} holds. Then $g$ is not Eisenstein.
\end{lemma}

\begin{proof}
    Suppose $g$ is Eisenstein. Then $g$ is a linear combination of $\mathbb 1$ and $(\frac{\cdot}{p}) \circ N$. We have that for all $q \nmid pN\ell$, \begin{align*}
        &T_q \left(\mathbb 1\right) = (q+1) \cdot \mathbb 1, &T_q \left(\left(\frac{\cdot}{p}\right) \circ N\right) = (q+1)\left(\frac{q}{p}\right) \cdot \left(\frac{\cdot}{p}\right)\circ N,
    \end{align*}
    so we must have that $\lambda_q(g) = \Tr \, \overline \rho_{\widetilde{f}}(\Frob_q) = (q+1)(\frac{q}{p})^e$ for $e\in \{0,1\}$. It follows by Chebotarev that $$\overline \rho_{\widetilde f}^{\rm ss} \cong \eta \oplus \eta\,\omega,$$
    where the superscript ${\rm ss}$ denotes semisimplification, $\eta$ is the Galois character corresponding to $(\frac{\cdot}{p})^e$ and
    $\omega$ is the Teichm\"uller character. This implies that $\overline \rho_{\widetilde f}$ is reducible, which violates \eqref{item_CR_2}.
\end{proof}

\begin{lemma} \label{lemma_Eisenstein}
    Let $q \neq  p$ be a prime at which $B'$ is split and let $t_q \in {B'_q}^\times$ be an element satisfying $N(t_q) =q$. Suppose that $h: {B'}^\times \backslash \widehat {B}'^\times \lra \mathcal O_1$ is invariant under right-translation by a compact open subgroup of $\widehat{B}'^\times$ and also under $\GL_2(\Z_q)$. For $\beta\in\mathcal O_1$, if $(1+\beta t_q)\cdot h$ is Eisenstein, then $h$ is Eisenstein.
\end{lemma}
\begin{proof}
    This is \cite[Lemma 5.4]{chidahseihanticycloLvaluescrelle}.
\end{proof}

Let ${\rm CM}: = \widehat \jmath(K^\times) \backslash \widehat B'^\times$ and let $\Gal(K^{\rm ab}/K)$ act on ${\rm CM}$ by ${\rec_K(a)}\cdot [b] = [\widehat \jmath(a)b]$ (\textit{cf}. \eqref{eqn_definite_Galois_action}), while $B'^\times_p$ acts by right multiplication.

\begin{defn}
    An element $\sigma \in \Gal(K^{\rm ab}/ K)$ is $p$-rational if $\sigma = \rec_K(a)$ for some $a \in \widehat K^\times$ satisfying $a_p \in K^\times \Q_p^\times \subset K_p^\times$.
\end{defn}

\begin{theorem}[Cornut-Vatsal]\label{theorem_Cornut-Vatsal}
    Let $U' \subset \widehat B'^\times$ be a compact open subgroup, $\mathcal{X} \subset {\rm CM}$ a $B'^\times_p$-orbit, and $R \subset \Gal(K^{\rm ab}/K)$ a finite set consisting of elements which are pairwise distinct modulo the subgroup of $p$-rational elements. Then for all but finitely many $x$ in the image of $\mathcal{X}$ in ${\rm CM}/U'$, we have
    \begin{equation*}
        \{(\tau r \cdot  x)_{r \in R} \in \prod_{r \in R} X(U'): \tau \in \Gal(K^{\rm ab}/K)\} = \{(y_r)_{r\in R} \in \prod_{r \in R}X(U'): \exists \tau, \, N(y_r) = N(\tau r \cdot x) \, \forall r \in R\}
    \end{equation*}
\end{theorem}

\begin{proof}
    This is \cite[Corollary 2.10]{cornutvatsal_2005} with $\mathfrak{S} = \emptyset$, $\mathfrak{R} = R$, $\mathcal{G} = \Gal(K^{\rm ab}/K)$, $H = U'$ and $P=p$. 
\end{proof}

\begin{proposition}
    $\theta(\mathcal{D'}) \neq 0$ in $\Lambda_1(\Gamma_\ac)$.
\end{proposition}

\begin{proof}
    By definition (\textit{cf.} \S \ref{subsubsec_theta_els}), it suffices for us to show that there exists some $m\geq 1$ and $\tau \in \Gal(H_m/K)$ such that
    \begin{equation}
        \sum_{\delta \in \Delta_{\ac}} g(\delta \tau \cdot x_m) \neq 0, \label{eqn_goal_nonzero}
    \end{equation}
    where we have set $x_m := x_m(1)$. We will proceed by way of contradiction, where the main input will be Theorem \ref{theorem_Cornut-Vatsal}.
        \item[\textit{Step 1:}] We start by laying out the set-up. Let $S$ denote the set of primes $q \neq p$ that are ramified in $K$, so $q\mathcal{O}_K = \mathfrak{q}^2$ for some prime ideal $\mathfrak{q}\subset \mathcal{O}_K$. Since $\Frob_\mathfrak{q}^2$ acts trivially on $H_\infty$, we have an inclusion $\Delta' \subset \Delta_\ac$, where $\Delta' \subset \Gal(H_\infty/K)$ is the subgroup generated by the $\Frob_{\mathfrak{q}}$ for such $q$. Let $S_0 \subset S$ be such that $\{\Frob_\mathfrak{q}\}_{\mathfrak{q} \in S_0}$ is an $\FF_2$-basis for $\Delta'$ and set 
        $$\tau_\fq := \widehat \jmath(\varpi_\q) \in B'^\times_q,$$
        for a fixed uniformizer $\varpi_\fq \in \mathcal{O}_{K_\fq}$. Since $\varsigma_{m,q}=1$ and $\jmath(K^\times_q)$ is commutative, we have that $\Frob_\q \cdot [\widetilde{\jmath}(a) \varsigma_m] = [\widetilde{\jmath}(a)\varsigma_m \tau_q]$ for all $a \in \widehat K^\times$, and hence, for every $y$ in the Galois orbit of $x_m$,
        \begin{equation*}
            \sum_{\delta' \in \Delta'}g(\delta' \cdot y) = \prod_{q \in S_0}\left((1+ \tau_q)\cdot g\right)(y)=: \vartheta(y),
        \end{equation*}
        where $\vartheta$ is invariant under right translation by the group $U_0' := U' \cap \bigcap_{q \in S_0} \tau_q U' \tau_q^{-1}$. Note that $N(U_0') = N(U')$ so the image of the induced map on $X(U'_0)$ is still $Z$.  If we let $R \subset \Gal(K^{\rm ab}/K)$ be a set of representatives for the elements of $\Delta_\ac / \Delta'$, then we can rewrite \eqref{eqn_goal_nonzero} as 
        \begin{equation}
            \sum_{r\in R} \vartheta(\tau  r\cdot x_m) \neq 0. \label{eqn_goal_rewrite}
        \end{equation}

        \item[\textit{Step 2:}] We now show that the elements of $R$ are pairwise distinct modulo $p$-rational elements. To that end, suppose that $\delta \in  \Delta_\ac$ is the image of some $p$-rational element $\rec_K(d)$. Since $d_p \in K^\times \Q^\times_p$, which acts trivially on $H_\infty$ under the reciprocity map, we may assume that $d_p =1$. Since $\delta^t =1$ for some $t$, we can write for each $m\geq 1$,  $d^t = k_m \cdot u_m $ for some $k_m \in K^\times$, $u_m \in \widehat {\mathcal O}_K^\times [p^m]$ with 
        $$k_m = u_{m,p}^{-1} \in (\Z_p + p^m \mathcal{O}_{K,p})^\times.$$
        Since $k_1 k_m^{-1} = u_mu_1^{-1} \in \mathcal{O}_K^\times$ for all $m$, we have that 
        $$k_1 \in \bigcap_m \mathcal{O}_K^\times \cdot (\Z_p + p^m\mathcal{O}_{K,p} ) = \mathcal{O}_K^\times \cdot \Z_p^\times,$$
        and hence, $k_1 \in \mathcal{O}_K^\times \cdot \Q^\times$. If $\mathfrak{d}$ denotes the ideal corresponding to $\delta$ by class field theory, then it follows that $\mathfrak{d}^t = (k_1)$ and $\overline{\mathfrak{d}} = \mathfrak{d}$, which can only occur if $\mathfrak{d}$ is a product of rational primes and primes that ramify in $K$. By definition, this implies that $\delta \in \Delta'$.
        \item[\textit{Step 3:}] $\vartheta$ is not Eisenstein. Indeed, if $\vartheta$ is Eisenstein, then so is $g$ by the definition of $\vartheta$ together with Lemma \ref{lemma_Eisenstein}, which contradicts Lemma \ref{lemma_g_not_Eis}.
        \item[\textit{Step 4:}] We now apply Theorem \ref{theorem_Cornut-Vatsal} with $\mathcal{X}$ equal to the $B_p'^\times$-orbit of $x_1$, so that $\{x_m\}_{m\geq 1}\subset \mathcal{X}$ and with $x= x_m$ for $m$ large enough such that the equality in \textit{loc. cit.} holds. 

        If $D_K \neq p$, then $(\frac{\cdot}{p}) \circ N$ is not the trivial character on $\widehat K^\times$, so the composition
        \begin{equation*}
            \widehat K^\times \lra X(U'_0) \xrightarrow{N} Z
        \end{equation*}
        is surjective, and the Galois orbit of $x_m$ intersects non-trivially both fibres of $N: X(U'_0) \lra Z$. 
        Now suppose \eqref{eqn_goal_rewrite} fails. By Theorem \ref{theorem_Cornut-Vatsal}, this is equivalent to saying that for all $\tau \in \Gal(H_m /K)$, and for all $(y_r)_{r\in R}$ with $N(y_r) = N(\tau r\cdot x_m)$,
        $$\sum_{r \in R}\vartheta(y_r) = 0. $$
        Varying $y_{r_0}$ for some fixed $r_0$ shows that $\vartheta$ is constant on the fibre $N^{-1}\left(N(\tau r_0 \cdot x_m)\right)$, and this implies that $\vartheta$ is Eisenstein by the discussion at the beginning of this paragraph, contradicting Step 3. 

        If $D_K=p$, then $S = \emptyset $, $\vartheta = g$, $U'_0 =U'$ and $(\frac{\cdot}{p}) \circ N$ is the trivial character. We again assume that \eqref{eqn_goal_rewrite} fails for all $\tau$. Arguing as in the previous case, this implies that $g$ is constant on $N^{-1}(1)$,  where now the Galois orbit of $x_m$ only intersects this fibre. By Remark~\ref{rk:cr2} and Chebotarev, we may find an inert prime $q\nmid Np\ell$ such that $\lambda_q(g) \neq  0$. Moreover, if $b \in X(U'_0)$ satisfies $N(b) =-1$, then 
        $$\left(\frac{\cdot}{p}\right) \circ N(b\pi_i) = (\frac{q}{p})N(b)= 1, \qquad i \in \{0, \dots, q\}, $$
        where $\iota(\pi_i) = \begin{pmatrix}
            q & i \\
            0 & 1
        \end{pmatrix}$ for $i \in \{0, \dots, q-1\}$ and $\iota(\pi_q) = \begin{pmatrix}
            1 & 0 \\
            0 &q
        \end{pmatrix}$, so that the $\pi_i$ are representatives of the double coset decomposition defining $T_q$. Hence we have the equality 
        $$g(b) = \lambda^{-1}_q(g)\cdot (T_qg)(b) = \lambda_q^{-1}(g)(q+1)c$$
        if $g$ is sent to $c$ on $N^{-1}(1)$, which is independent of $b$, and thus $g$ must be Eisenstein, contradicting Lemma~\ref{lemma_g_not_Eis}. 
\end{proof}

\subsection{Kolyvagin's conjecture}\label{s:kolyvagin}
Kolyvagin's conjecture from \cite{kolyvagin91} concerns the non-vanishing of the system of derived Heegner classes. The classes in question are of the form
\[
c_n\in H^1(K,T_{\widetilde f}/I_nT_{\widetilde f}),
\]
where \(n\) is a square-free product of Kolyvagin primes and
\[
I_n=\sum_{\ell\mid n}I_\ell,
\qquad
I_\ell=(\ell+1,a_\ell(\widetilde f)).
\]
The conjecture asserts that this system is non-trivial. It was proved by W.~Zhang \cite{zhang14} in the good ordinary case under certain residual ramification hypotheses. More recent results, under various hypotheses, have been obtained in \cite{SkinnerZhang,BCGS,kim24,kim26,sweeting,LPV,DaRonche,CastellaSano,leizhao}. With the exception of \cite{SkinnerZhang}, which treats the case of multiplicative reduction at \(p\), these results assume that \(p\) is a prime of good reduction. 

We explain the relation between Kolyvagin's conjecture and the equality criterion discussed in Theorem~\ref{them_indef_main_conjecture}(iii).
Following \cite[Lemma~3.6 and Proposition 3.7]{BCK}, consider the following residual non-vanishing condition on the bipartite Euler system:
\[
\theta(\mathcal D_{\mathfrak n}) 
\not\equiv 0
\mod{(\varpi,\gamma_{\mathrm{ac}}-1)},
\]
for some $\mathfrak{n}\in\mathfrak{N}_1^\mathrm{def}$.
In other words, at least one theta element has nonzero image in
$\Lambda(\Gamma_{\mathrm{ac}})/(\varpi,\gamma_\ac-1)$.
In \cite{BCK}, the authors showed that in the good ordinary case, this condition implies Howard's criterion for the equality of the Heegner point main conjecture. They make use of Zhang's proof of Kolyvagin's conjecture, where a specific definite level-raised form satisfying this condition is constructed. The following proposition gives the corresponding implication in our setting.

\begin{proposition}
    Assume \eqref{item_Heeg}, \eqref{item_ord}, \eqref{item_CR} and \eqref{item_irr} hold. Suppose that for every $n\ge1$, the image of the system given by Lemma~\ref{lemma_Bipartite_Euler_System_Exists} under $\Lambda(\Gamma_\ac)\to \cO_1$ (the evaluation map at $\mathbb{1}$ composed with the projection $\cO\to\cO_1$)  gives rise to a free bipartite Euler system of odd type $\ES_n$ for $(T_1,\mathfrak{F},\mathfrak{L}_n)$.
Let $\m_{\Lambda(\Gamma_\ac)}=(\varpi,\gamma_\ac-1)$ denote the maximal ideal of $\Lambda(\Gamma_\ac)$. If
\[
 {\theta(\mathcal D_{\mathfrak n_0}^{(1)})} \not\equiv0\mod{\mathfrak m_{\Lambda(\Gamma_\ac)}}
 \qquad\text{for some }\mathfrak n_0\in\mathfrak N_1^{\mathrm{def}},
\]
then the equality criteria in Theorems~\ref{thm_definite_main_conj}(iii) and Theorem~\ref{them_indef_main_conjecture}(iii) hold with $k_0=1$ for every height-one prime $\fP$ of $\Lambda(\Gamma_\ac)$.
\end{proposition}
\begin{proof}
By Lemma~\ref{lemma_Bipartite_Euler_System_Exists}, evaluation at the trivial character followed by reduction modulo $\varpi$ gives
\[
 \ES_n=\ES_1|_{\fN_n},  \qquad n\geq1.
\]

    For $\fa\in\fN_1$, put
    \[
    r(\fa)=\dim_{\cO_1}\Sel_{\Delta_\fa}(K,T_1).
    \]
    Fix $n\ge1$ and choose $\fa\in\fN_n$ such that $r(\fa)$ is minimal. If $r(\fa)>0$, it follows from Theorem~\ref{theorem_Howard_hyp_e} and \cite[Corollary~2.2.10]{howard06} that there exists $\ell\in\fL_n$ with $\ell\nmid\fa$, $a\ell\in\fN_n$ such that
    \[
    r(\fa\ell)=r(\fa)-1,
    \]
    which contradicts the minimality of $r(\fa)$. Thus, $r(\fa)=0$. In particular, $\fa\in\fN_n^\mathrm{def}$.

    Our hypotheses imply that $\ES_1$ is non-zero. As $\cO_1$ is a field, we can apply \cite[Theorem~2.5.1]{howard06} with $\delta=0$, giving
    \[
    {\theta(\mathcal{D}_\fa^{(1)})}\not\equiv0\mod \fM_{\Lambda(\Gamma_\ac)}.
    \]
By the compatibility with $\ES_n$,
\[
 {\theta(\mathcal{D}_\fa^{(n)})}  \equiv{\theta(\mathcal{D}_\fa^{(1)})} \not\equiv0\mod{\fM_{\Lambda(\Gamma_{\ac})}}.
\]   
As $(\fP,\varpi)\subseteq\fM_{\Lambda(\Gamma_\ac)}$ for any height-one prime, this implies Howard's equality criteria.
\end{proof}

In \cite{BCK}, the non-vanishing of ${\theta(\mathcal{D}_{\fn_0}^{(1)})}$ is obtained from Zhang's proof of Kolyvagin's conjecture in the good ordinary case. Here we write ${\theta(\mathcal{D}_\fn^{(j)})}$ for the elements denoted by $\lambda_j(\fn)$ in \cite{BCK}.
Using a common coefficient ring, put
\[
 r(h):=\dim_{\mathcal O_1}
       \operatorname{Sel}(K,A_h[\varpi]),
 \qquad A_h=V_h/T_h.
\]
Under the hypotheses of \cite[Proposition~3.7]{BCK}, the level-raising argument in the proof of
\cite[Theorem~9.1]{zhang14} produces congruent forms $f_0,\ldots,f_s$ satisfying
\[
 r(f_{i+1})=r(f_i)-2,
 \qquad r(f_s)=1.
\]
The proof of \cite[Theorem~7.2]{zhang14} then gives a further level-raised form $g$, corresponding to a definite admissible product $\fn_0$, with
\[
 \operatorname{Sel}(K,A_g[\varpi])=0,
 \qquad
 \operatorname{Sel}(K,A_g)=0.
\]
Normalizing $L^{\mathrm{alg}}(g/K)$ by Gross's period, \cite[Theorems~7.1 and~6.4]{zhang14} give
\[
 \operatorname{ord}_{\varpi}L^{\mathrm{alg}}(g/K)
 =
 \operatorname{length}_{\mathcal O}
       \operatorname{Sel}(K,A_g)
 =0.
\]
Let $\mathbf 1$ denote the trivial character of $\Gamma_{\mathrm{ac}}$.
Evaluating at $\mathbf 1$, the interpolation calculation in \cite[Proposition~3.7]{BCK} gives
\[
 \mathbf 1\bigl({\theta(\mathcal{D}_{\fn_0}^{(1)})}\bigr)^2
 =
 \overline{u\,e_p(g,\mathbf 1)L^{\mathrm{alg}}(g/K)}
 \neq0
 \quad\text{in }\mathcal O_1,
\]
where $u\in\mathcal O^\times$ and $e_p(g,\mathbf 1)\in\mathcal O^\times$ by the nonanomalous
hypothesis. Hence, we obtain
\[
 {\theta(\mathcal{D}_{\fn_0}^{(1)})}\not\equiv0\mod{\mathfrak m_{\Lambda(\Gamma_\ac)}}.
\]
We hope that future work on Kolyvagin's conjecture in the potentially ordinary setting will establish the non-vanishing condition above, which would imply that the divisibilities in Theorems~\ref{thm_definite_main_conj} and \ref{them_indef_main_conjecture} are equalities.

\subsection{Arithmetic applications}\label{s:app}
Throughout, we assume that \eqref{item_Heeg}, \eqref{item_ord}, \eqref{item_CR} and \eqref{item_irr} hold. In what follows, we write $\Sel(K,M)$, $M\in\{T_{\widetilde{f}},V_{\widetilde{f}},A_{\widetilde{f}}\}$ for the Bloch--Kato Selmer group. The local Bloch--Kato subgroups are denoted by $H^1_\f(K_v,M)$.
\begin{lemma}\label{lem:BK}
   For $v|p$, we have $H^1_\f(K_v,A_{\widetilde{f}})=H^1_\ord(K_v,A_{\widetilde{f}})$. Moreover, there is an exact sequence
   \[
   0\to\Sel(K,A_{\widetilde{f}})\to\Sel_{N^-}(K,A_{\widetilde{f}})\to\bigoplus_{v|N^+} H^1_{\mathrm{fin}}(K_v,A_{\widetilde{f}}),
   \]
   where the last term is finite.
\end{lemma}
\begin{proof}
    For $v|p$, Lemma~\ref{lemma_vanishing_cohomology_at_p}, \cite[Corollary~3.8.4]{blochkato}, and the local Euler characteristic formula give 
    \[
    H^1_\f(K_v,V_{\widetilde{f}})=\image\left(H^1(K_v,F_p^+V_{\widetilde{f}})\lra H^1(K_v,V_{\widetilde{f}})\right).
    \]
    Since $H^2(K_v,F_p^+T_{\widetilde{f}})=0$ by Lemma~\ref{lemma_vanishing_cohomology_at_p}, we have 
    \[
    H^1(K_v,F_p^+V_{\widetilde{f}})\twoheadrightarrow H^1(K_v,F_p^+A_{\widetilde{f}}). 
    \]
    Taking images in  $H^1(K_v,A_{\widetilde{f}})$, we have
    \[
H^1_\f(K_v,A_{\widetilde{f}})=\image\left(H^1(K_v,F_p^+A_{\widetilde{f}})\to H^1(K_v,A_{\widetilde{f}})\right)=H^1_\ord(K_v,A_{\widetilde{f}}),
    \]
    which proves the first assertion. Away from $p$, the local comparison is the same as in \cite[proof of Corollary~6.15]{chidahsiehanticyclomainconjformodformscomposito}: the local conditions defining the two Selmer groups differ only at $v\mid N^+$, where they are $0$ and $H^1_\mathrm{fin}(K_v,A_{\widetilde{f}})$, respectively. The proof that the latter group is finite is the same as \textit{loc. cit.} This implies the second assertion of the lemma.  
\end{proof}

For $v|N^+$, set $t_v=\length_\cO H^1_\mathrm{fin}(K_v,A_{\widetilde{f}})$. We have the following application of Theorem~\ref{thm_definite_main_conj} in the definite case.

\begin{theorem}\label{thm:app1}
    Suppose $\omega(N^-)$ is odd. If $L(\widetilde{f}/K,1)\ne 0$, then the Bloch--Kato Selmer group $\Sel(K,A_{\widetilde{f}})$ is finite. Furthermore,
    \[
    \length_\cO\Sel(K,A_{\widetilde{f}})+\sum_{v|N^+}t_v\le\ord_\varpi\mathbb{1}(\mathscr{L}_f).
    \]
\end{theorem}
\begin{proof}
Take $g=f$,  $\psi=\epsilon^{-1}$ and $\chi=\mathbb{1}$ in \S\ref{sec:KP}. If $\epsilon^2\ne\mathbb{1}$, or $\epsilon^2=\mathbb{1}$ and $p\nmid D_K$, then $\widetilde{\mathbb{1}}=\epsilon^{-1}\circ\mathrm{N}_{K/\QQ}$ is ramified above $p$, so the first interpolation formula in Theorem~\ref{theorem_p-adic_l-function_interpolation} applies.
If $\epsilon^2=\mathbb{1}$ and $p|D_K$, the second interpolation formula applies. The factor $1-a_p(f)^{-1}\widetilde\chi(\fp)$ is a $p$-adic unit under \eqref{item_CR_1}(ii) since $\widetilde\chi(\fp)\in\{\pm1\}$. Furthermore,  $L_p(\pi_{f^\circ},\widetilde{\mathbf1},1)$
is finite and nonzero. Hence, in both cases,
%Then $\widetilde\chi$ is ramified at every prime above $p$ by the inertia argument in the proof of Lemma~\ref{lemma_vanishing_cohomology_at_p}(i) \red{maybe we need to assume $p \nmid D_K$ when $\epsilon^2 = \mathbb 1$?}. In particular, $\chi$ is good. It follows from Theorem~\ref{theorem_p-adic_l-function_interpolation} that
$$ L(\widetilde{f}/K,1)\ne 0\quad \Longleftrightarrow\quad \mathbb{1}(\mathscr{L}_f)\ne0.$$ 
Thus, by Proposition~\ref{prop_comparing_thetas} and the definition of $\mathscr{L}_f$, we have
\[
\mathbb{1}(\mathscr{L}_f)=u\theta_0^2\ne0,
\]
where $\theta_0=\mathbb{1}(\theta_\infty)$ and $u$ is a $p$-adic unit.    
    
    Combining the exact sequence in Lemma~\ref{lem:BK} with the Poitou--Tate exact sequence gives  \begin{equation}
0\to\Sel(K,A_{\widetilde{f}})\to\Sel_{N^-}(K,A_{\widetilde{f}})\to\bigoplus_{v|N^+} H^1_{\mathrm{fin}}(K_v,A_{\widetilde{f}})\to \Sel(K,T_{\widetilde{f}})^\vee,\label{eq:ses}
      \end{equation} 
     since $H^1_\f(K_v,T_{\widetilde{f}})$ and $H^1_\f(K_v,A_{\widetilde{f}})$ are orthogonal complements of each other under the local Tate pairing for every finite place $v$ by \cite[Proposition~3.8]{blochkato}.  
Applying Proposition~\ref{prop_lengths_of_Selmer}(i) with $\fP=(\gamma_\ac-1)$ and $k$ sufficiently large implies that $\Sel_{N^-}(K,A_{\widetilde f})$ is finite, and
    \[\length_\cO\Sel_{N^-}(K,A_{\widetilde f})\le
    \length_\cO\Sel_{N^-}(K,A_{\widetilde f})+2\delta_\fP(k)=2\ord_\varpi(\theta_0)=\ord_\varpi\mathbb{1}(\mathscr{L}_f).
    \]

As $\Sel(K,A_{\widetilde f})$ is finite, $\Sel(K,T_{\widetilde{f}})$ has rank zero. As $H^0(K,T_{\widetilde{f}}/\varpi T_{\widetilde{f}})=0$, we have $\Sel(K,T_{\widetilde{f}})=0$. The exact sequence \eqref{eq:ses} gives
\[
\length_\cO\Sel_{N^-}(K,A_{\widetilde{f}})=\length_\cO\Sel(K,A_{\widetilde{f}})+\sum_{v|N^+}t_v.
\]
Thus, the desired inequality follows.
\end{proof}
    
 We now consider the indefinite case. Set $\kappa_0=\mathbb{1}(\kappa_\infty)\in\Sel_{N^-}(K,T_{\widetilde{f}})$.   

\begin{theorem}\label{thm:app2}
   Suppose that $\omega(N^-)$ is even. Assume that $\kappa_0\ne 0$. Then, the $\cO$-module $\Sel_{N^-}(K,T_{\widetilde{f}})$ is free of rank one, whereas $\Sel(K,A_{\widetilde{f}})$ is of corank one. Furthermore, the quotient $C_K:=\Sel_{N^-}(K,A_{\widetilde{f}})/\Sel(K,A_{\widetilde{f}})$ is finite, and
    \[
    \length_\cO\left(\Sel(K,A_{\widetilde{f}})^\vee_{\mathrm{tor}}\right)+\length_\cO(C_K)\le 2\length_\cO\left(\Sel_{N^-}(K,T_{\widetilde{f}})/\cO\kappa_0\right).
    \]
\end{theorem}
\begin{proof}
The first assertion follows from applying Proposition \ref{prop_lengths_of_Selmer}(ii) with $\fP=(\gamma_\ac-1)$, and $k$ chosen so that $\kappa_0\notin\varpi^k\Sel_{N^-}(K,T_{\widetilde{f}})$. 
Furthermore,
\[
\length_\cO\left(\Sel_{N^-}(K,A_{\widetilde{f}})^\vee_{\mathrm{tor}}\right)\le\length_\cO\left(\Sel_{N^-}(K,A_{\widetilde{f}})^\vee_{\mathrm{tor}}\right)+2\delta_\fP(k)= 2\length_\cO\left(\Sel_{N^-}(K,T_{\widetilde{f}})/\cO\kappa_0\right).
\]

As $H^1_{\mathrm{fin}}(K_v,A_{\widetilde{f}})$ is finite for all $v|N^+$ by Lemma~\ref{lem:BK}, the quotient $C_K$ is also finite. In particular, we have the short exact sequence 
$$
0\to C_K^\vee\to\Sel_{N^-}(K,A_{\widetilde{f}})^\vee_\mathrm{tor}\to \Sel(K,A_{\widetilde{f}})^\vee_\mathrm{tor}\to 0.
$$
Hence, the final assertion follows.
\end{proof}

For an elliptic curve $E/\mathbf Q$ as in the introduction, let $c_v(E/K)$ denote the local Tamagawa number at $v$.
If $v\mid\ell\mid N^+$, then $K_v=\QQ_\ell$.
Since $p>3$, Tate's algorithm shows that $p\mid c_v(E/K)$ can occur only at a split multiplicative prime with $p\mid\ord_\ell(\Delta_E)$, which is excluded by \eqref{item_CR_3}(ii).
Consequently, the local groups in the exact sequence of Lemma~\ref{lem:BK} vanish, so $C_K=0$. Theorems~\ref{thm:app1} and \ref{thm:app2} therefore give:

\begin{corollary}The following properties hold.
\begin{itemize}
    \item[i)] Suppose that $\omega(N^-)$ is odd and $L(E/K,1)\ne0$. Then $\rank_\ZZ E(K)=0$ and $\Sha(E/K)[p^\infty]$ is finite. Furthermore, 
    \[
    \ord_p\#\Sha(E/K)[p^\infty]\le \ord_p(\mathbb{1}(\mathscr{L}_f)).
    \]
    \item[ii)] Suppose that $\omega(N^-)$ is even. Assume that the image of $\kappa_0$ in $H^1(K,V_p(E))$ is a non-zero element of $E(K)\otimes\Qp\hookrightarrow H^1(K,V_p(E))$. Then $\rank_\ZZ E(K)=1$ and $\Sha(E/K)[p^\infty]$ is finite. Furthermore,
    \[
    \ord_p\#\Sha(E/K)[p^\infty]\le 2\ord_p\#(\Sel_{N^-}(K,T_p(E))/\Zp\kappa_0).
    \]
\end{itemize}
    \end{corollary}

\bibliographystyle{amsalpha}
\bibliography{references}
\end{document}